\documentclass{amsart}
\usepackage{amsthm,latexsym,amsfonts,amsmath,amssymb,mathrsfs,array,manfnt,stmaryrd}
\usepackage[pdftex]{graphicx}
\usepackage[all]{xy}
\usepackage{paralist}
\usepackage{framed}
\usepackage{cancel}
\usepackage{enumitem}
\usepackage{asymptote}
\usepackage{units}
\usepackage{setspace}
\usepackage{amscd}
\usepackage{microtype}
\usepackage{tikz}
\usetikzlibrary{matrix}
\usetikzlibrary{calc}
\usepackage{tikz-cd}

\newcommand{\into}{\hookrightarrow}

\newcommand{\pderiv}[3][]{\frac{\partial^{#1} #2}{\partial {#3}^{#1}}}

\newcommand{\abs}[1]{\left\lvert#1\right\rvert}

\newcommand{\Z}{\ensuremath{\mathbb{Z}}}

\newcommand{\C}{\ensuremath{\mathbb{C}}}
\newcommand{\Q}{\ensuremath{\mathbb{Q}}}

\newcommand{\A}{\ensuremath{\mathbb{A}}}

\newcommand{\M}{\mathcal{M}}
\newcommand{\Mbar}{\overline{M}}
\newcommand{\Hbar}{\overline{\mathcal{H}}}

\newcommand{\isom}{\cong}

\renewcommand{\mod}{\thickspace\mathrm{mod}\thinspace}
\renewcommand{\P}{\ensuremath{\mathbb{P}}}
\renewcommand{\H}{\ensuremath{\mathcal{H}}}

\newcommand{\floor}[1]{\left\lfloor#1\right\rfloor}

\renewcommand{\bar}[1]{\overline{#1}}

\DeclareMathOperator{\Aut}{Aut}

\DeclareMathOperator{\Spec}{Spec}

\DeclareMathOperator{\Per}{Per}

\DeclareMathOperator{\val}{val}

\DeclareMathOperator{\codim}{codim}
\DeclareMathOperator{\st}{st}

\usepackage{color}

\newcommand{\margincolor}{red}      
\definecolor{darkgreen}{rgb}{0,0.7,0}

\newcounter{margincounter}
\newcommand{\marginnum}{
\ifnum\value{margincounter}<10
\textcolor{\margincolor}{\begin{picture}(0,0)\put(2.2,2.4){\circle{9}}\end{picture}\footnotesize\arabic{margincounter}}
\else\ifnum\value{margincounter}<100
\textcolor{\margincolor}{\begin{picture}(0,0)\put(4.256,2.5){\circle{11}}\end{picture}\footnotesize\arabic{margincounter}}
\else
\textcolor{\margincolor}{\begin{picture}(0,0)\put(6.8,2.5){\circle{14}}\end{picture}\footnotesize\arabic{margincounter}}
\fi\fi
}

\newcommand{\cC}{\mathcal{C}}

\newcommand{\cE}{\mathcal{E}}

\newcommand{\cP}{\mathcal{P}}

\newcommand{\bA}{\mathbf{A}}
\newcommand{\bB}{\mathbf{B}}

\newcommand{\bV}{\mathbf{V}}

\theoremstyle{plain}
\newtheorem{theorem}{Theorem}
\numberwithin{theorem}{section}

\newtheorem{thm}[theorem]{Theorem}
\newtheorem{algorithm}[theorem]{Algorithm}
\newtheorem*{thm*}{Theorem}

\newtheorem{cor}[theorem]{Corollary}

\theoremstyle{definition}

\newtheorem{caution}[theorem]{Caution}

\theoremstyle{remark}

\newtheorem{rem}[theorem]{Remark}
\newtheorem{obs}[theorem]{Observation}

\usepackage{tikz}
\usetikzlibrary{matrix}
\usepackage{fullpage}
\usepackage{bbm}

\DeclareMathOperator{\CR}{CR}
\usepackage{hyperref}

\newcommand{\cm}{\mathrm{cm}}

\DeclareMathOperator{\fl}{fl}
\DeclareMathOperator{\mk}{mk}

\begin{document}
\title{Equations at infinity for critical-orbit-relation families of
  rational maps}
\date{\today}
\author{Rohini Ramadas and Rob Silversmith}

\begin{abstract}
  We develop techniques for using compactifications of Hurwitz spaces to study families of rational maps $\P^1\to\P^1$ defined by critical orbit relations. We apply these techniques in two settings: We show that the parameter space $\Per_{d,n}$ of degree-$d$ bicritical maps with a marked 4-periodic critical point is a $d^2$-punctured Riemann surface of genus $\frac{(d-1)(d-2)}{2}$. We also show that the parameter space $\Per_{2,5}$ of degree-2 rational maps with a marked 5-periodic critical point is a 10-punctured elliptic curve, and we identify its isomorphism class over $\Q$. We carry out an experimental study of the interaction between dynamically defined points of $\Per_{2,5}$ (such as PCF points or punctures) and the group structure of the underlying elliptic curve. 
\end{abstract}
\maketitle

\section{Introduction}\label{sec:Intro}

\subsection{COR varieties} The $1$-parameter family of quadratic polynomials $\{f_c(z)=z^2+c\}_{c\in\C}$ is the most-studied family of complex dynamical systems; it is the home of the Mandelbrot set. The parameter space $\C$ (with coordinate $c$) can be thought of as the space parametrizing, up to conjugacy, all quadratic rational self-maps of the Riemann sphere with a distinguished critical fixed point (the point at $\infty$). Milnor \cite{Milnor1993}, Stimson \cite{Stimson1993}, Rees \cite{Rees2009}, and others (elaborated below) have defined and studied natural generalizations, including the $1$-parameter families $\Per_{2,n}$ of quadratic rational maps with a (distinguished) critical point that is periodic of period exactly $n$. More generally, if one fixes a degree $d$, and imposes $r$ ``independent" conditions on orbits of critical points, one obtains a $(2d-2-r)$-parameter family of degree-$d$ rational functions with marked critical points whose orbits are as specified; these are called \emph{critical-orbit-relation (COR) families} or (when referring to the parameter space of such a family) \emph{COR varieties}. COR varieties are affine algebraic varieties defined over $\Q$; fundamental transversality results from Teichm\"uller theory, due to Thurston and others \cite{DouadyHubbard1993,Epstein2009}, imply that they are smooth (with a well-understood class of exceptions).

COR varieties are fundamental in the study of complex and arithmetic dynamics. A conjecture of Baker and DeMarco, proved in some cases \cite{BakerDeMarco2013,FavreGauthier2015, DeMarcoWangYe2015}, states that (again with certain exceptions) COR varieties are the only families of rational maps with a Zariski-dense subset of \emph{post-critically finite/(PCF)} maps, i.e. maps for which every critical point is (pre-)periodic. Also, \emph{deformation spaces} of rational maps, defined by Epstein \cite{Epstein2009}, arise naturally from the perspective of Teichm\"uller theory and complex dynamics and are covering spaces of COR families \cite{Rees2009,HironakaKoch2017,Hironaka2019,FirsovaKahnSelinger2016}.

COR varieties have been most studied in the $1$-parameter case $r=2d-3$. These \emph{COR curves} are punctured Riemann surfaces/smooth algebraic curves defined over $\Q$; one can therefore study their irreducible components, genera, number of punctures, gonalities\footnote{Recall that the \emph{gonality} of a Riemann surface $X$ is the minimum degree of a nonconstant holomorphic map $X\to\P^1$}, isomorphism classes over $\Q$, and so on --- and how these invariants vary in natural countable collections of COR curves like $\Per_{d,n}$. Here are four cases where questions of this type have been addressed --- the last is the focus of this paper:
\begin{itemize}
    \item The COR curves $X_n$ parametrizing cubic polynomials with an $n$-periodic critical point. By works of Arfeux, Bonifant, DeMarco, Kiwi, Milnor, Pilgrim, and Schiff \cite{BonifantKiwiMilnor2010,DeMarcoPilgrim2011,DeMarcoSchiff2012,ArfeuxKiwi2020}, $X_n$ is known to be irreducible, with algorithms given for its genus and number of punctures.
    \item The ``dynatomic'' curves $Y_{d,k,\ell}$ parametrizing maps $f:z\mapsto z^d+c$ with a marked preperiodic point $z_0$ of type $(k,\ell)$, so $f^k(z_0)=f^{k+\ell}(z_0)$). (These are technically not COR curves but a mild generalization.) Doyle-Poonen \cite{DoylePoonen2020} have shown that the curves $Y_{d,k,\ell}$ are irreducible, and that there are finitely many of them with any fixed gonality.
    \item The COR curves $Z_{2,k},Z_{3,k}$ parametrizing quadratic rational maps or cubic polynomials with a preperiodic critical point of type $(k,1)$. Buff-Epstein-Koch \cite{BuffEpsteinKoch2018} have shown that $Z_{2,k}$ and $Z_{3,k}$ are irreducible.
    \item The COR curves $\Per_{d,n}$ parametrizing degree-$d$ \emph{bicritical} maps with a marked $n$-periodic critical point. These are comparatively mysterious --- it is known that for $n\le 4$, $\Per_{2,n}$ is rational, so admits a global parametrization (which can be written explicitly), while $\Per_{2,5}$ has genus $1$ \cite{Stimson1993}.
\end{itemize}
COR varieties naturally embed in the moduli spaces $\M_d$ of degree-$d$ rational maps, and are often studied via this perspective (see Section \ref{sec:whyhurwitz}).
Alternatively, COR varieties also embed naturally into Hurwitz spaces (\cite{Epstein2009}), again see Section \ref{sec:whyhurwitz}); this was used by Hironaka-Koch \cite{HironakaKoch2017} and Hironaka \cite{Hironaka2019} to study the topology of the deformation space of $\Per_{2,4}$. In addition, Harris-Mumford \cite{HarrisMumford1982} have constructed \emph{admissible covers} compactifications of Hurwitz spaces; these compactifications have a richly combinatorial boundary stratification and local coordinates. 

\subsection{Summary of results} In this work, we develop machinery that uses the geometry and combinatorics of spaces of admissible covers to study COR varieties by producing explicit defining equations ``at infinity''. We do this explicitly for the COR curves $\Per_{d,n}$ (\textbf{Algorithm \ref{alg:Punctures}}), but our methods can be adapted for general COR varieties (see Remark \ref{rem:Generalize}). We apply the methods to two cases: $\Per_{d,4}$ for varying $d$, and $\Per_{2,5}$ (which, as mentioned, has genus 1). First, we describe their punctures explicitly as $1$-parameter families of rational maps degenerating to maps between trees of $\P^1$ (see Figures \ref{fig:StrataThatIntersectC5} and \ref{fig:TableOfBoundaryStrataBicritical}, also compare to \cite{DeMarcoPilgrim2011}). We use this to show: 

\begin{thm}\label{thm:GenusFormula0}
 $\Per_{d,4}$ is isomorphic to a smooth plane curve of degree $d$, punctured at $d^2$ points. In particular, $\Per_{d,4}$ has genus $\frac{(d-1)(d-2)}{2}$ and gonality $d-1$.
\end{thm}

Note that the curves $\Per_{d,4}$ are ``special'', in the sense that a general smooth curve of genus $\frac{(d-1)(d-2)}{2}$ does not admit an embedding as a smooth plane curve. As a result, the gonality of $\Per_{d,4}$ is much lower than the gonality $\floor{\frac{d^2-3d+8}{4}}$ of a general curve of genus $\frac{(d-1)(d-2)}{2}$. It would be interesting to know if there is a general sense in which COR curves are ``special'' points of their moduli spaces $M_{g}$.

We next show
\begin{thm}\label{thm:whatisc5}
$\Per_{2,5}$ is isomorphic to the $\Q$-elliptic curve $\cC_{17a4}$, punctured at 10 points.
\end{thm}
The elliptic curve identifier 17a4 is from the \emph{L-functions and modular forms database (LMFDB)}, which lists various invariants of the curve, e.g. that its conductor and discriminant are equal to $17$, its $j$-invariant is equal to $\frac{3^3 11^3}{17}$, its Mordell-Weil group is isomorphic to $\Z/4\Z$, and its endomorphism ring is isomorphic to $\Z$. The curve $\cC_{17a4}$ therefore contains both points of number-theoretic interest, e.g. its countable set of torsion points, and points of dynamical interest, e.g. its finite set of punctures, and its countable set of PCF points. These are all defined over $\bar{\Q}$. We find:
\begin{thm}\label{thm:c5torsion}
The $10$ punctures of $\Per_{2,5}$ include the 4 $\Q$-rational points of $\cC_{17a4}$. The other $6$ punctures are points of infinite order on $\cC_{17a4}$. The $20$ PCF points of $\Per_{2,5}$ that parametrize maps with both critical points in the same periodic $5$-cycle are points of infinite order on $\cC_{17a4}$.
\end{thm}

As far as we are aware, this is the first such investigation into a dynamical moduli space that is also an elliptic curve. A recurring question in arithmetic dynamics is to classify which conditions on critical points are realizable over $\Q$, see \cite{LukasManesYap2014,Poonen1998,BrezinByrneLevyPilgrimPlummer2000}. We have the following immediate consequence of Theorem \ref{thm:c5torsion}:
\begin{cor}
No quadratic rational function $\P^1\to\P^1$ with $\Q$-coefficients has a 5-periodic critical point.
\end{cor}

We have also drawn a tiled analog of the Mandelbrot set inside $\Per_{2,5}$, via the identification of $\widetilde{\Per}_{2,5}$ with a parallelogram in the upper half-plane of $\C$ (see Section \ref{sec:Mandelbrot}). The result is Figure \ref{fig:Mandelbrot} below. Punctures of $\Per_{2,5}$ are labelled in red; observe that they all appear to be at ``cusps'' of the fractal. Analogous fractals have been depicted for $\Per_{2,2}$ by Timorin \cite{Timorin2006}, for $\Per_{2,3}$ by Rees in \cite{Rees2009}, and for $\Per_{2,4}$ by Gage-Jackson in \cite{GageJackson2011}. 

\begin{figure}
    \centering
    \includegraphics[width=.9\textwidth,height=1.014\textwidth]{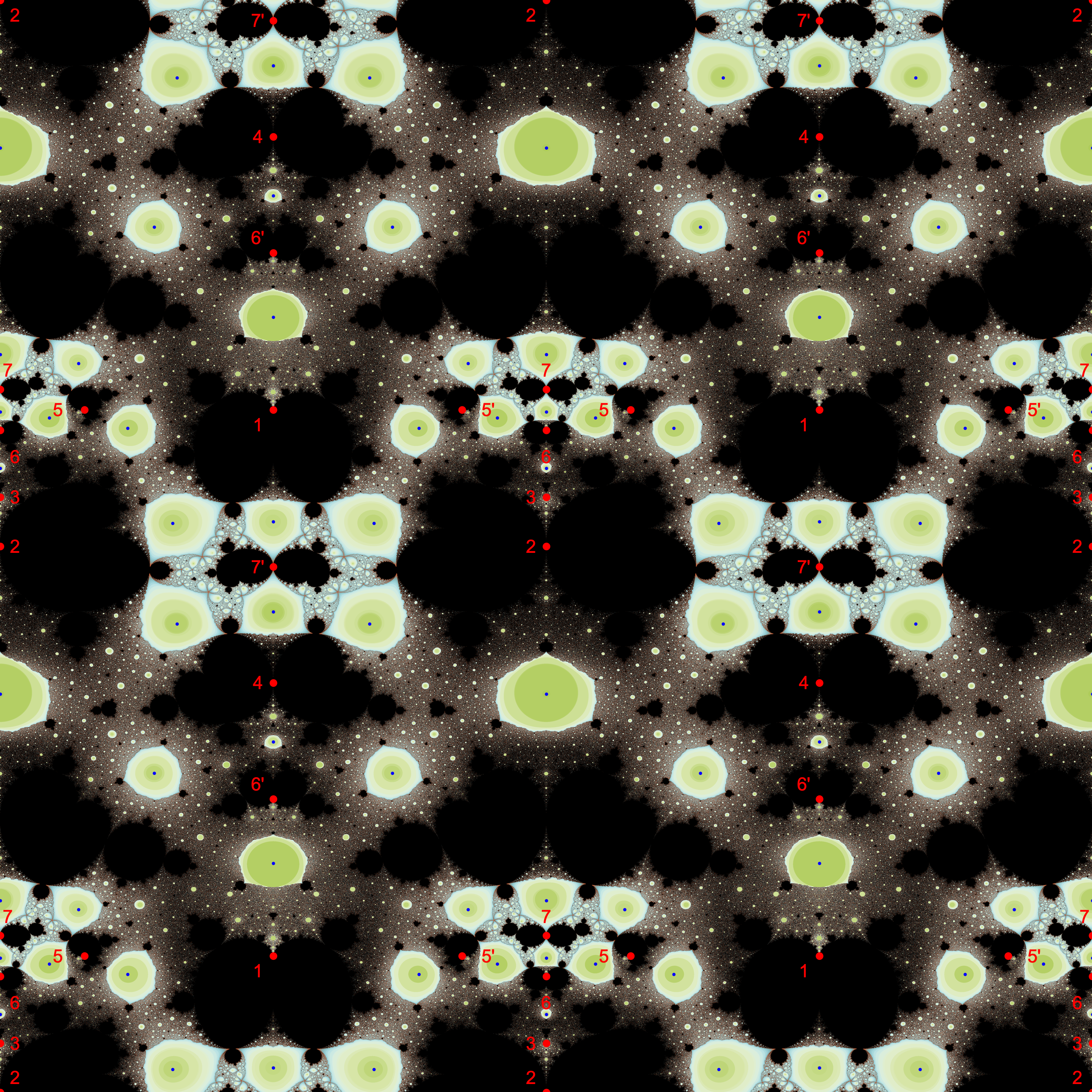}
    \caption{An analog of the Mandelbrot set in $\Per_{2,5}$. Points are colored black if they correspond to maps for which the second critical point is not in the attracting basin of the critical 5-cycle. Punctures are labeled in red, with numbering as per Table \ref{fig:StrataThatIntersectC5} and Section \ref{sec:C5InHBar}. Small blue dots in the centers of some of the colored regions denote PCF maps for which both critical points are in the same 5-cycle.}
    \label{fig:Mandelbrot}
\end{figure}

\subsection{Comparison of ambient moduli spaces}\label{sec:whyhurwitz}
As noted above, $\Per_{d,n}$ admits several natural maps that are helpful in studying its geometry. We discuss some aspects, advantages and disadvantages of each.

\subsubsection*{The space of rational maps} The moduli space $\M_d$ of conjugacy classes of degree-$d$ rational maps on $\P^1$ is a $(2d-2)$-dimensional affine variety defined over $\Q$ \cite{Silverman1998}. (Silverman \cite{Silverman1998} and Schmitt \cite{Schmitt2017} have constructed GIT compactifications of $\M_d$.) The space $\M_d$ admits a finite ramified cover $\M_d^{\cm}$ parametrizing degree-$d$ rational maps with all critical points marked. Any COR variety parametrizing degree-$d$ rational functions is naturally embedded in $\M_d^{\cm}$ as a codimension-$r$ subvariety, where $r$ is the number of constrained critical points\footnote{For precise statements on the structure of these constraints, called ``postcritical portraits'', see \cite[Sec. 9]{FloydParryPilgrim2017}; we do not need these to study $\Per_{d,n}$.}. In particular we have a natural inclusion $\Per_{d,n}\into\M_d^{\cm}$, and the resulting composition $\Per_{d,n}\to\M_d$, which is injective except at points $(f,p)\in\Per_{d,n}$ where \emph{both} critical points are $n$-periodic; these are singular points of the image. Unfortunately, these embeddings have very high degree, as well as bad behavior near punctures. For example, $\M_2\cong\A^2$, and the curve $\Per_{2,5}$ (which has genus 1 and has 10 nodes in $\M_2$, see Theorem \ref{thm:whatisc5}) is defined by a degree-15 equation \cite{Milnor1993}. In the compactification $\P^2\supseteq\M_2,$ this means the closure $\overline{\Per_{2,5}}$ must have singularities along the line at infinity that are equivalent to 80 simple nodes (as a smooth degree-15 curve in $\P^2$ has genus 91). Since $\Per_{2,5}$ has only 10 punctures (again by Theorem \ref{thm:whatisc5}), these singularities must be quite nasty. Note however that Stimson \cite{Stimson1993} has proved that the branches of these singularities are smooth.

\subsubsection*{The space of bicritical maps} $\Per_{d,n}$ lies inside the subvariety $\M_d^{\mathrm{bicrit}}\subseteq\M_d$ parametrizing bicritical rational maps of $\P^1$, i.e. those with exactly two critical points. (Of course, $\M_2^{\mathrm{bicrit}}=\M_2$.) Milnor \cite{Milnor2000} showed that $\M_d^{\mathrm{bicrit}}\cong\C^2$, and therefore studied $\Per_{d,n}$ (the normalization of) a plane curve (singular, as noted just above). Milnor considers the compactification $\P^2$ of $\C^2,$ and studies the punctures of $\Per_{d,n}$ by intersecting $\Per_{d,n}$ with lines nearer and nearer to the line at $\infty$.
\subsubsection*{The space of maps with a marked periodic point} Manes \cite{Manes2009} introduced the moduli space $\M_d(n)$ of degree-$d$ rational maps of $\P^1$ with a marked $n$-periodic point; $\M_d(n)$ is a finite ramified cover of $\M_d$. Blanc-Canci-Elkies \cite{BlancCanciElkies2013} studied $\M_2(n)$, showing that it is rational if $n\le5$. A pair $(f,p)\in\M_2(n)$ has an associated \emph{multiplier} $\lambda\in\C$, defined by $(f^n)'(p)$. This gives a fibration $\M_2(n)\to\C,$ and $\Per_{2,n}$ is precisely the fiber over $0\in\C$. In particular, it follows from \cite{BlancCanciElkies2013} that $\M_2(5)$ is naturally isomorphic to $M_{0,5}$, where $M_{0,n}$ denotes the moduli space of $n$ marked curves on $\P^1$ up to change of coordinates. The resulting embedding $\Per_{2,5}\into M_{0,5}$ is the same as the one we consider in Section \ref{sec:CubicCurve}.
\subsubsection*{Moduli spaces of point-configurations on $\P^1$} More generally, any COR subvariety admits maps to $M_{0,n}$ by marking $n$ points in the constrained critical orbit. In particular, $\Per_{d,n}$ admits a natural map to $M_{0,n}$ sending a rational map $f$ to the (ordered) configuration of $n$ points in the periodic critical orbit. In turn, $M_{0,n}$ admits a smooth projective 
compactification $\Mbar_{0,n}$ whose \textit{boundary} $\Mbar_{0,n}\setminus M_{0,n}$ parametrizes degenerate point-configurations, i.e. configurations of marked points on nodal trees of $\P^1$s. The local geometry of $\Mbar_{0,n}$ at the boundary is described via combinatorics, see Section \ref{sec:M0nBarCoords}. Degenerations within a COR variety (such as punctures of $\Per_{d,n}$), therefore, have corresponding degenerations of point-configurations.

\subsubsection*{Hurwitz spaces}
In order to understand the above degenerations of point-configurations, it is helpful is to keep track of all of critical points and critical values, not just those in the constrained critical orbits. The natural way to do this is via Hurwitz spaces, which parametrize maps of curves with prescribed branching behavior and marked critical points. The map from a COR variety to $M_{0,n}$ factors through a Hurwitz space $\H$ (Section \ref{sec:Setup}); while the map to $M_{0,n}$ is not necessarily an embedding, the map to $\H$ is an embedding \cite{Epstein2009}. Degenerations within COR varieties occur when critical values (postcritical points) collide with each other. Correspondingly, the Harris-Mumford ``admissible covers'' compactification $\Hbar$ of $\H$ exactly describes the degenerations of maps that occur when critical values collide (Section \ref{sec:ModuliSpaces}). Importantly, Harris and Mumford give a combinatorial description of the deformation theory of such degenerations --- i.e. of the local geometry of $\Hbar$ (Section \ref{sec:HnBarCoords}) --- and this description is crucial in allowing us to give equations near punctures. These key features --- that $\Hbar$ captures degenerations of configurations of critical values, and that its local geometry can be expressed combinatorially in terms of these degenerations --- seems to be unique to the admissible covers compactification. As a result, it is the authors' opinion that moduli spaces of admissible covers are an important and natural tool for studying COR varieties.


While the perspective of admissible covers is well-suited to the study of punctures, one unfortunate aspect is that COR varieties have very high codimension in Hurwitz spaces and in $M_{0,n}$, when the number of marked points is large. This can make the global geometry of COR varieties challenging to study; in our examples, the codimension of the embedding is small enough that we are nonetheless able to make conclusions about global geometry: $\Per_{d,4}$ is a hypersurface in a Hurwitz space, and $\Per_{d,5}$ is a hypersurface in $M_{0,5}$.



We make an additional note: It is natural to study a COR \emph{curve} via mapping to $\Mbar_{0,n}$, then forgetting all but 4 marked points to obtain a map to $\Mbar_{0,4}\cong\P^1$. This collection of maps to $\P^1$ produces bounds on the genus and gonality of the curve, and admissible covers allow one to compute their degrees purely combinatorially. In particular, from Sections \ref{sec:C5InHBar} and \ref{sec:GenusFormula}, we have:
\begin{cor}
 The minimal-degree map $\Per_{d,n}\to\P^1$ is achieved by a forgetful map to $M_{0,4}$ in the cases $(d,n)=(2,5)$ and $(d,n)=(*,4)$.
\end{cor}

\subsubsection*{Parametrizations} In a few small cases, one may study $\Per_{d,n}$ not via an embedding but via an explicit parametrization, see e.g. \cite{Timorin2006} for $\Per_{2,2}$, \cite{Rees2009} for $\Per_{2,3}$, \cite{GageJackson2011} for $\Per_{2,4}$. This is done by constructing a normal form for such a rational function, e.g. $f(z)=\frac{z+c}{z^2}$ for $(f,0)\in\Per_{2,2}$, where $c\ne0$. This is already impossible for $\Per_{2,5},$ which is not rational. (This limitation was observed in \cite{Stimson1993}.)

\subsection*{Acknowledgements}
Both authors thank Sarah Koch and Joseph H. Silverman for useful conversations; the first author also thanks Xavier Buff, Laura DeMarco, and Rob Benedetto.

\subsection{Funding}
The first author was supported by an NSF postdoctoral fellowship (DMS-1703308), and a Tamarkin Assistant Professorship at Brown University. The second author was supported by an RTG-Zelevinsky Research Instructorship at Northeastern University, part of NSF grant DMS-1645877.

\subsection{Notation and conventions}\label{sec:Conventions} We work over a field $k$ of characteristic zero. If $f_1,\ldots, f_n$ are regular functions on a scheme $U$, we denote by $Z(f_1,\ldots f_n)$ the subscheme of $U$ obtained as their common vanishing locus. If $f_1$ and $f_2$ are regular functions on a scheme $U$, then $Z(f_1)=Z(f_2)$ if and only if there is a regular function $\alpha$, non-vanishing on $U$, such that $f_1=\alpha f_2$; in this case we write $f_1\sim_{U} f_2$, or $f_1\sim f_2$ if $U$ is clear from context.

We write $[n]:=\{1,\ldots,n\}.$ For a tree $\sigma,$ we write $\mathbf{V}(\sigma)$ and $\mathbf{E}(\sigma)$ for its vertex and edge sets, respectively. If $\mathbf{S}$ is a finite set, an \textit{$\mathbf{S}$-marked tree} is a tree $\sigma$ together with a map $\mk:\mathbf{S}\to\mathbf{V}(\sigma).$ The elements of $\mathbf{S}$ are called \textit{legs} or \textit{markings}. If $\mk(p)=v$, we say $p$ is incident to $v$; the leg $p$ should be visualized as a half-edge incident to $v$. A \textit{flag of $\sigma$ at $v$} is an element $e\in\mathbf{E}(\sigma)\cup S$ with the property that $e$ is incident to $v$. For $v\in\mathbf{V}(\sigma),$ we denote by $\mathbf{F}(v)$ the set of flags of $\sigma$ at $v$. We define the \textit{valence} $\val(v)=\abs{\mathbf{F}(v)}$. Given a vertex $v$ of $\sigma$, and $p\in \mathbf{S}$, the \textit{flag connecting $v$ to $p$} is defined to be $p$ if $p$ is incident to $v$, and is otherwise defined to be the unique edge $e$ incident $v$ with the property that $v$ and $p$ are in distinct connected components of $\sigma\setminus\{e\}$.


\section{Moduli spaces of marked curves and maps}\label{sec:ModuliSpaces}
We introduce several moduli spaces, all defined over $\Q$, and
consider the relationships between them. 

\subsection{Setup}\label{sec:Setup} 
Let $d\ge2.$ We define $\Per_{d,n}$ to be the moduli space parametrizing pairs $(f,p),$ where $f:\P^1\to\P^1$ is a degree-$d$ rational function with exactly two critical points, $p\in\P^1$ is one of the critical points, and $p$ is a periodic point of $f$ with period exactly $n$, up to conjugation by $\Aut(\P^1)$. It follows from standard results \cite{DouadyHubbard1993,Epstein2009,Silverman1998} that $\Per_{d,n}$ is a smooth affine curve defined over $\Q$.

\begin{caution}
The notation $\Per_{d,n}$ should not be confused with various similarly-notated objects in complex dynamics, e.g. $\Per_{d,n}(\lambda)$ for the family of degree-$d$ rational functions with an $n$-periodic critical point of multiplier $\lambda$.
\end{caution}

In this paper, we study $\Per_{d,n}$ via a natural rational map to a Hurwitz space, which we now describe. Recall that for $n\ge 3$, the moduli space $M_{0,n}$ of $[n]$-marked (or $n$-marked) rational curves is a smooth $(n-3)$-dimensional variety over $\Q$ that parametrizes tuples $(C, p_1,\ldots p_n)$, where $C$ is a smooth projective genus-zero curve and $p_1,\ldots, p_n\in C$ are distinct points. Note that since $C$ has a rational point it is isomorphic to $\P^1$. 

For fixed $n\ge2$, we fix sets $$\mathbf{A}_n=\{a_*,a_1\}\sqcup\{a_{i,k}\}_{\substack{2\le i\le n\\0\le k\le d-1}},\quad\quad\mathbf{A}_{n,0}=\{a_*,a_1\}\sqcup\{a_{i,0}\}_{2\le i\le n},\quad\quad\text{and}\quad\quad\mathbf{B}_n=\{b_1,\ldots,b_n,b_*\}.$$ (So $\mathbf{A}_{n,0}\subseteq\mathbf{A}_n.$ We define a map $\phi:\mathbf{A}_n\to\mathbf{B}_n$ by $\phi(a_*)=b_*,$ $\phi(a_1)=b_2,$ and $\phi(a_{i,k})=b_{i+1\mod n}$ for $i=2,\ldots,n$ and $k=0,\ldots,d-1$. We also define a degree map $\deg:\mathbf{A}_n\to\{1,2\}$ by $\deg(a_*)=\deg(a_1)=d,$ and $\deg(a_{k,i})=1$ for $i=2,\ldots,n$ and $k=0,\ldots,d-1$.
(This notation will be convenient in Sections \ref{sec:HnBar} and \ref{sec:HnBarCoords}.) We define the Hurwitz space $\H_{d,n}$ to be the
moduli space of
tuples $(C,D,f),$ where
\begin{itemize}
\item $C$ is a smooth genus zero curve, together with an injective map $\iota:\bA_{n,0}\into C$. We suppress the notation of $\iota$ and write $a\in C$ for the image, under $\iota$, of $a\in \bA_{n,0}$.
\item $D$ is a smooth genus zero curve, similarly marked by the set $\bB_n$, and
\item $f$ is a degree-$d$ map $f:C\to D$ such that for all $a\in\mathbf{A}_{n,0}$, $f(a)=\phi(a)$ and $f$ has local degree $\deg(a)$ at $a$ (and $f$ has no other ramification).
\end{itemize}
Note that $\H_{d,n}$ parametrizes maps up to independent changes of coordinates on $C$ and $D$, thus the behavior, under iteration, of $f\in \H_{d,n}$ is not well-defined. The space $\H_{d,n}$ is a smooth quasiprojective $(n-2)$-dimensional variety over $\Q$ \cite{Fulton1969,RomagnyWewers2006}. It admits morphisms $\pi_1,\pi_2:\H_{d,n}\to M_{0,n}$, where $\pi_1(C,D,f)=(D,b_1,\ldots,b_{n})$, and $\pi_2(C,D,f)=(C,a_1,a_{2,0},\ldots,a_{n,0})$.

Let $\Delta_n\subseteq M_{0,n}\times M_{0,n}$ denote the diagonal. Given a point of $(\pi_1\times\pi_2)^{-1}(\Delta_n)$, there is a unique identification of $C$ with $D$ that identifies $a_1$ with $b_1$ and $a_{i,0}$ with $b_i$ for $i=2,\ldots, n$. Under this identification of source and target, $f$ is a dynamical quadratic rational map whose critical point $a_1$ is $n$-periodic, so we have a map $(\pi_1\times\pi_2)^{-1}(\Delta_n)\to\Per_{d,n}.$ This map has a clear inverse on the locus $\Per_{d,n}^{\circ}\subseteq\Per_{d,n}$ where the unmarked critical point is not in the orbit of the marked critical point. Thus $(\pi_1\times\pi_2)^{-1}(\Delta_n)\cong\Per_{d,n}^{\circ}$.
(Note that $\Per_{d,n}^{\circ}$ is dense in $\Per_{d,n}$, again by \cite{DouadyHubbard1993,Epstein2009,Silverman1998}.)

\begin{rem}\label{rem:BaseChange}
    It will be important to consider the base change of $\H_{d,n}$ to $\Q(\zeta),$ where $\zeta$ is a primitive $d$th root of unity. After this base change, we may canonically mark additional points $a_{i,k}\in\mathbf{A}_n$ on the source curve $C$, as follows. For fixed $2\le i\le n,$ we may identify $C$ with $\P^1$ so that $a_*=\infty$, $a_1=0,$ and $a_{i,0}=1$, and identify $D$ with $\P^1$ so that $b_*=\infty$, $b_2=0,$ and $b_{i+1\mod n}=1$. With respect to these coordinates, $f:C\to D$ is identified with the function $z\mapsto z^d$. We then mark each additional preimage $\zeta^k\in\P^1$ by the element $a_{i,k}\in\mathbb{A}_n;$ this marking is canonical over $\Q(\zeta),$ though not over $\Q$ if $d>2.$
    
    These extra markings define a factorization of $\pi_2:\H_{d,n}\to M_{0,n}$ through $M_{0,\bA_n}$; this will allow us to define additional coordinate functions on $\H_{d,n}$ and its compactification $\Hbar_{d,n}$, see Sections \ref{sec:HdnCoords} and \ref{sec:HnBarCoords}.
\end{rem}

\subsection{The Deligne-Mumford compactification}\label{sec:M0nBar} By works of Knudsen, Deligne-Mumford, and Grothendieck, there is a smooth projective compactification $\Mbar_{0,n}$ of $M_{0,n}$ parametrizing \textit{stable $n$-marked genus-zero curves} \cite{Knudsen1983,DeligneMumford1969}. A stable $n$-marked genus-zero curve is a connected nodal curve $C$ of arithmetic genus zero, together with distinct marked smooth points $p_1,\ldots,p_n$ of $C$, such that every irreducible component of $C$ contains at least three ``special points'' -- a special point is a node or a marked point.

To any stable $n$-marked genus-zero curve $(C,p_1,\ldots,p_n)$ is associated its $n$-marked dual graph $\sigma$ (see \cite{}). The fact that $C$ is stable of genus zero implies that $\sigma$ is a \textit{stable $n$-marked tree}, i.e. a tree each of whose vertices has valence at least 3. (Recall from Section \ref{sec:Conventions} that the valence of a vertex includes marked legs.) The assignment $C\mapsto\sigma$ defines a stratification on $\Mbar_{0,n}$ by locally closed strata; these strata are in bijection with stable $n$-marked trees. We denote by $Q_\sigma$ the locally closed stratum $\{C\thickspace|\thickspace C \text{ has dual tree }\sigma\}$, and we denote by $\bar{Q_{\sigma}}$ its Zariski closure. A curve $C'$ is in $\bar{Q_{\sigma}}$ if and only if there exists a set of edges of its dual tree $\sigma'$ such that contracting those edges yields $\sigma.$

By regarding each irreducible component of a stable curve $C$ as a stable curve in its own right, we have natural isomorphisms 
\begin{align}\label{eq:QSigmaProduct}
    Q_\sigma\cong\prod_{v\in\mathbf{V}(\sigma)}M_{0,\val(v)}\quad\quad\text{ and }\quad\quad\bar{Q_\sigma}\cong\prod_{v\in\mathbf{V}(\sigma)}\Mbar_{0,\val(v)}.
\end{align} We therefore have $\dim(Q_\sigma)=\sum_{v\in\mathbf{V}(\sigma)}(\val(v)-3)$. If $\sigma$ is a stable tree and $e$ is an edge of $\sigma$,  or equivalently $\codim(Q_\sigma)=\abs{\mathbf{E}(\sigma)}.$

Let $(C,p_1,\ldots,p_n)$ be a (possibly unstable) curve with distinct marked points, and let $\sigma$ be its $n$-marked dual tree. There is a \textit{stabilization} $\st(C,p_1,\ldots,p_n)\in\Mbar_{0,n}$ obtained by contracting irreducible components with $\le2$ special points until a stable curve is reached. Correspondingly, the stabilization of $\st(\sigma)$ of $\sigma$ is a stable $n$-marked tree obtained by repeatedly choosing a vertex $v$ with $\val(v)<3$, choosing an edge $e$ incident to $v$, and contracting $e$; this process results in the dual tree of $\st(C,p_1,\ldots,p_n)$. The process of stabilization induces \textit{forgetful morphisms} $\Mbar_{0,n}\to\Mbar_{0,n-1}$; these extend the usual forgetful morphisms $M_{0,n}\to M_{0,n-1}.$ Composing them, one may forget any subset of $\{p_1,\ldots,p_n\}$ of size at most $n-3$.

\subsection{Compactifications of Hurwitz spaces}\label{sec:HnBar} Harris and Mumford \cite{HarrisMumford1982} constructed compactifications of Hurwitz spaces by moduli spaces of \textit{admissible covers}. These spaces parametrize maps of possibly nodal curves. We denote by $\Hbar_{d,n}$ the admissible covers compactification of $\H_{d,n}$; it parametrizes tuples $(C,D,f)$, where
\begin{itemize}
\item $C$ is stable $\bA_{n,0}$-marked genus zero curve,
\item $D$ is a stable $\bB_n$-marked genus zero curve.
\item $f$ is a finite degree-$d$ map $f:C\to D$ such that 
\begin{itemize}
\item For all $a\in\bA_{n,0}$, $f(a)=\phi(a)$ and $f$ has local degree $\deg(a)$ at $a$ (and $f$ has no other ramification at smooth points),
\item nodes of $C$ map to nodes of $D$ and smooth points of $C$ map to smooth points of $D$, and
\item (\textit{balancing condition}) at each node $\eta\in C$, the two branches of the node $\eta$ map to the two branches of the node $f(\eta)\in D$ with equal local degree, see \cite{HarrisMumford1982}.
\end{itemize}
\end{itemize}

By \cite{HarrisMumford1982}, $\Hbar_{d,n}$ is a smooth projective variety containing $\H_{d,n}$ as a dense open subset. (Generally, moduli spaces of admissible covers may be singular orbifolds --- though the singularities can be ignored in a certain sense, see \cite{AbramovichCortiVistoli2003} --- but the local coordinates given in \cite{HarrisMumford1982} imply that $\Hbar_{d,n}$ is smooth and is a variety.) There are natural maps $\pi_1,\pi_2:\Hbar_{d,n}\to\Mbar_{0,n}$ that extend the maps $\pi_1,\pi_2:\H_{d,n}\to M_{0,n}$, defined by: $\pi_1(C,D,f)=\st(D,b_1,\ldots,b_n)$, and $\pi_2(C,D,f)=\st(C,a_{1,0},\ldots,a_{n,0})$. 

\begin{rem}\label{rem:BaseChange2}
It is straightforward to check that the extra markings of Remark \ref{rem:BaseChange} extend to $\Hbar_{d,n}$ (again after base change to $\Q(\zeta)$). Again, this defines a factorization of $\pi_2$ through $\Mbar_{0,\bA_n}$.
\end{rem}

Given $(C,D,f)\in\Hbar_{d,n}$, we can extract its \textit{combinatorial type} $\gamma=(\sigma, \tau,\phi,\deg),$
where
\begin{itemize}
    \item $\sigma$ is the $\bA_{n,0}$-marked dual tree of $C$,
    \item $\tau$ is the $\bB_n$-marked dual tree of $D$,
    \item $\phi:\sigma\to\tau$ is a graph homomorphism; that is, a pair $(\phi_1,\phi_2),$ where $\phi_1:\mathbf{V}(\sigma)\to\mathbf{V}(\tau)$ and
    $\phi_2:\mathbf{E}(\sigma)\to\mathbf{E}(\tau)$ if $e\in\mathbf{E}(\sigma)$ is incident to $v\in\mathbf{V}(\sigma)$, then $\phi_2(e)\in\mathbf{E}(\tau)$ is incident to $\phi_1(v)\in\mathbf{V}(\tau),$ and
    \item $\deg:\mathbf{V}(\sigma)\cup\mathbf{E}(\sigma)\to\{1,d\}$ encodes, for each irreducible component of $C$, the degree of $f$ on that component, and for each node of $C$, the local degree of $f$ on each branch of the node.
\end{itemize}
For convenience, we denote both $\phi_1$ and $\phi_2$ by $\phi$. (Note that in Section \ref{sec:Setup} we defined maps $\phi:\mathbf{A}_{n,0}\to\mathbf{B}_n$ and $\deg:\mathbf{A}_{n,0}\to\{1,d\}$; this notation allows us to treat legs and edges uniformly.)
The definition of an admissible cover implies the following properties of a combinatorial type:
\begin{itemize}
    \item The graph map $\phi$ is surjective on both edges and vertices.
    \item If $a\in\mathbf{A}_{n,0}$ is a marking on the vertex $v\in\mathbf{V}(\sigma),$ then $\phi(a)\in\mathbf{B}_n$ is a marking on $\phi(v)\in\mathbf{V}(\tau).$
    \item For any $e\in\mathbf{E}(\tau)$, $\phi^{-1}(e)$ consists of two edges of $\sigma$, counted with multiplicity $\deg$.
    \item For any $w\in\mathbf{V}(\tau)$, $\phi^{-1}(w)$ consists of two vertices of $\sigma$, counted with multiplicity $\deg$.
    \item The degree-2 vertices $v\in\mathbf{V}(\sigma)$ are exactly those on the unique minimal path from the vertex marked by $b_*$ to the vertex marked by $b_2$.
    \item If $v\in\mathbf{V}(\sigma)$ has degree 2, then exactly two flags incident to $v$ have degree 2.
\end{itemize}

\begin{figure}
    \centering
    \begin{tikzpicture}
      \draw (-1.5,0) node {$\tau$};
\draw (-1.5,2) node
      {$\sigma$};
\draw (0,0)--(2,0);
\draw (0,0) node
      {$\bullet$};
\draw (0,0)--++(160:.4);
\draw (0,0)++(160:.4)
      node[left] {2};
\draw (0,0)--++(200:.4);
\draw (0,0)++(200:.4)
      node[left] {$*$};
\draw (1,0) node
      {$\bullet$};
\draw (1,0)--++(-110:.4);
\draw (1,0)++(-110:.4)
      node[below] {3};
\draw (1,0)--++(-70:.4);
\draw (1,0)++(-70:.4)
      node[below] {4};
\draw (2,0) node
      {$\bullet$};
\draw (2,0)--++(20:.4);
\draw (2,0)++(20:.4)
      node[right] {1};
\draw (2,0)--++(-20:.4);
\draw (2,0)++(-20:.4)
      node[right] {5};
\draw[->] (1,.75)--(1,0.25);
\draw
      (2,2.5)--(1,2.5)--(0,2)--(1,1.5)--(2,1.5);
\draw (2,2.5) node
      {$\bullet$};
\draw (2,2.5)--++(20:.4);
\draw (2,2.5)++(20:.4)
      node[right] {5};
\draw (2,2.5)--++(-20:.4);
\draw (2,2.5)++(-20:.4)
      node[right] {4};
\draw (1,2.5) node
      {$\bullet$};
\draw (0,2) node
      {$\bullet$};
\draw[very thick] (0,2)--++(160:.4);
\draw
      (0,2)++(160:.4) node[left] {1};
\draw[very thick]
      (0,2)--++(200:.4);
\draw (0,2)++(200:.4) node[left]
      {$*$};
\draw (1,1.5) node
      {$\bullet$};
\draw (1,1.5)--++(-110:.4);
\draw
      (1,1.5)++(-110:.4) node[below] {2};
\draw (1,1.5)--++(-70:.4);
      \draw (1,1.5)++(-70:.4) node[below] {3};
\draw (2,1.5) node
      {$\bullet$};
    \end{tikzpicture}
    \quad\quad\quad\quad\quad\quad
    \raisebox{12pt}{
    \begin{tikzpicture}
\draw (-1.5,2) node {$\bar{\sigma}$};
\draw (0,2)--(2,2);
\draw (0,2) node
      {$\bullet$};
\draw (0,2)--++(160:.4);
\draw (0,2)++(160:.4)
      node[left] {5};
\draw (0,2)--++(200:.4);
\draw (0,2)++(200:.4)
      node[left] {4};
\draw (1,2) node
      {$\bullet$};
\draw (1,2)--++(-90:.4);
\draw (1,2)++(-90:.4)
      node[below] {1};
\draw (2,2) node
      {$\bullet$};
\draw (2,2)--++(20:.4);
\draw (2,2)++(20:.4)
      node[right] {2};
\draw (2,2)--++(-20:.4);
\draw (2,2)++(-20:.4)
      node[right] {3};
    \draw (-1.5,0) node {$\bar{\tau}$};
\draw (0.5,0)--(1.5,0);
\draw (0.5,0) node
      {$\bullet$};
\draw (0.5,0)--++(140:.4);
\draw (0.5,0)++(140:.4)
      node[left] {2};
\draw (0.5,0)--++(180:.4);
\draw (0.5,0)++(180:.4)
      node[left] {3};
\draw (0.5,0)--++(220:.4);
\draw (0.5,0)++(220:.4)
      node[left] {4};
\draw (1.5,0) node
      {$\bullet$};
\draw (1.5,0)--++(20:.4);
\draw (1.5,0)++(20:.4)
      node[right] {1};
\draw (1.5,0)--++(-20:.4);
\draw (1.5,0)++(-20:.4)
      node[right] {5};
    \end{tikzpicture}
    }
    \caption{A boundary stratum in $\Hbar_{2,5}$ that does not intersect $\tilde{\Delta}_{2,5}$}
    \label{fig:NoIntersection}
\end{figure}
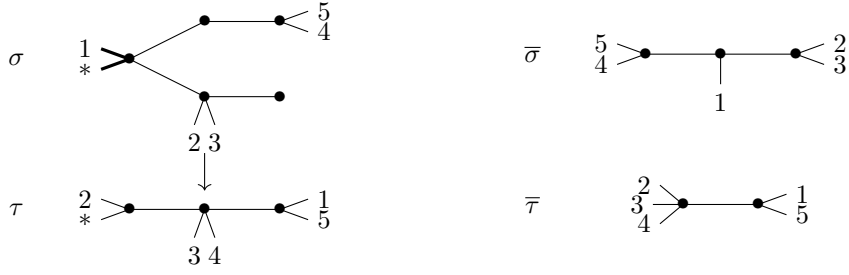

We visualize a combinatorial type by drawing $\phi$ as a vertical map; see Figure \ref{fig:NoIntersection} and Tables \ref{fig:StrataThatIntersectC5} and \ref{fig:DoNotIntersectC5}. For readability, 
we write simply $*$ or $i$ instead of $a_*,$ $a_1,$ $a_{i,0}$, $b_*,$ $b_i$. 

The assignment $(C,D,f)\mapsto\gamma$ defines a stratification on $\Hbar_{d,n}$ by locally closed strata. We denote by $R_\gamma$ the locally closed stratum corresponding to a combinatorial type $\gamma$; we refer to the Zariski closure $\bar{R_{\gamma}}$ as a closed stratum of $\Hbar_{d,n}$. As with $Q_\sigma$ and $\bar{Q_\sigma}$ above, we have decompositions
\begin{align}\label{eq:RGammaProduct}
    R_{\gamma}&\cong \prod_{\substack{w\in\bV(\tau)\\\abs{\phi^{-1}(w)}=1}}\H_{d,\val(w)-1}\times \prod_{\substack{w\in\bV(\tau)\\ \abs{\phi^{-1}(w)}=d}}M_{0,\val(w)}&\quad\quad\bar{R_{\gamma}}&\cong \prod_{\substack{w\in\bV(\tau)\\\abs{\phi^{-1}(w)}=1}}\Hbar_{d,\val(w)-1}\times \prod_{\substack{w\in\bV(\tau)\\ \abs{\phi^{-1}(w)}=d}}\Mbar_{0,\val(w)}.
\end{align}

We introduce special names $\bar{\sigma}$ and $\bar{\tau}$ for $\st(\sigma,a_{1,0},\ldots,a_{n,0})$ and $\st(\tau,b_1,\ldots,b_n)$ respectively; note that $\pi_1(R_\gamma)= Q_{\bar{\tau}}\subseteq\Mbar_{0,\bA_{n,0}}$ and $\pi_2(R_\gamma)\subseteq Q_{\bar{\sigma}}\subseteq\Mbar_{0,\bB_n}$.

Let $\bar{\Delta_n}$ denote the diagonal in $\Mbar_{0,n}\times \Mbar_{0,n}$, and let $\tilde{\Delta}_{d,n}:=(\pi_1\times \pi_2)^{-1}(\bar{\Delta_n})\subset \Hbar_{d,n}$. Then $\tilde{\Delta}_{d,n}$ is compact, and from above we have $\tilde{\Delta}_{d,n}\cap\H_{d,n}=\Per_{d,n}^{\circ}$. However, we will see that $\tilde{\Delta}_{d,n}$ is not, in general, equal to the Zariski closure of $\Per_{d,n}^{\circ}$; that is, $\tilde{\Delta}_{d,n}$ may have irreducible components contained in the boundary of $\Hbar_{d,n}$. We let $\overline{\Per_{d,n}}$ denote the Zariski closure of $\Per_{d,n}^{\circ}$ in $\Hbar_{d,n}$, so we have $\Per_{d,n}^{\circ}\subset\overline{\Per_{d,n}}\subset \tilde{\Delta}_{d,n}.$ We also let $\widetilde{\Per}_{d,n}$ be the normalization of $\overline{\Per_{d,n}}$, i.e. the unique smooth projective completion of $\Per_{d,n}$. Since $\Per_{d,n}^{\circ}$ is a curve, the complements  $\overline{\Per_{d,n}}\setminus\Per_{d,n}^{\circ}$ and $\widetilde{\Per}_{d,n}\setminus\Per_{d,n}^{\circ}$ are finite sets of points (over $\bar{\Q}$), with a surjective map $\widetilde{\Per}_{d,n}\setminus\Per_{d,n}^{\circ}\to \overline{\Per_{d,n}}\setminus\Per_{d,n}^{\circ}$.

\begin{obs}\label{obs:CombinatorialCondition}
  Let $(C,D,f)\in\tilde{\Delta}_{d,n}$ have combinatorial type $(\sigma,\tau,\phi,\deg).$ By definition, we have  $$\st(C,a_{1,0},\ldots,a_{n,0})=\st(D,b_1,\ldots,b_n)$$ as points of $\Mbar_{0,n},$ so in particular $\bar{\sigma}=\bar{\tau}$. This provides a combinatorial necessary condition for a boundary stratum $R_\gamma$ to intersect $\tilde{\Delta}_{d,n}.$ For example, Figure \ref{fig:NoIntersection} shows a combinatorial type with $n=5$ such that $\bar{\sigma}\ne\bar{\tau}$.
\end{obs}

\section{Cross-ratios and coordinates on moduli spaces}\label{sec:Coordinates} The goal of this section is introduce the tools we need for analyzing $\tilde{\Delta}_{d,n}\cap R_\gamma$ for various boundary strata $R_\gamma$ of $\Hbar_{d,n}$. For this purpose, it is sufficient to work in an infinitesimal neighborhood\footnote{Indeed, we \textit{will} need to use power series in our analysis; see e.g. Section \ref{sec:R6}, which involves computing Taylor series expansions of the square root function.} of $R_\gamma,$ i.e. a scheme locally of the form $\Spec(k[s_1,\ldots,s_q][[s_{q+1},\ldots,s_r]]),$ where $s_1,\ldots,s_q$ are local coordinates on $R_\gamma$, and $s_{q+1},\ldots,s_r$ are regular functions on $\Hbar_{d,n}$ that restrict to a basis of sections of the conormal bundle of $R_\gamma\subseteq\Hbar_{d,n}.$ We now describe such coordinates and how to use them.

\subsection{Cross-ratios and global coordinates on \texorpdfstring{$M_{0,n}$}{M\_\{0,n\}}}\label{sec:M0nCoords}
Given $(C, p_1,\ldots, p_4)\in M_{0,4}$, there is a unique $\lambda\in\P^1\setminus\{0,1,\infty\}$ such that $(C,p_1,\ldots,p_4)=(\P^1,\infty,0,1,\lambda)\in M_{0,4}$; the value $\lambda$ is the \textit{cross-ratio} of the four points $p_1,p_2,p_3,p_4$. The assignment $(C, p_1,\ldots, p_4\mapsto\lambda$ induces an identification $\CR(1,2,3,4):M_{0,4}\xrightarrow{\isom}\P^1\setminus\{0,1,\infty\}=\A^1\setminus\{0,1\}$, which extends to an identification $\Mbar_{0,4}\xrightarrow{\isom}\P^1$. Any 4-tuple of distinct elements $\{i_1,i_2,i_3,i_4\}\subseteq [n]$ determines a forgetful map to $\Mbar_{0,4}$, and, via the identification above, defines a cross-ratio map $\CR(i_1,\ldots, i_4):\Mbar_{0,n}\to\P^1$. When restricted to $M_{0,n}$, these maps take values in $\A^1\setminus\{0,1\}$, and satisfy relations arising from changing coordinates on $\P^1$. We refer to these as the \textit{cross-ratio functional equations}. For example,
\begin{align*}
    \CR(1,3,4,2)&=\frac{1}{1-\CR(1,2,3,4)}&&\text{and}&
    \CR(1,3,4,5)&=\frac{\CR(1,2,3,5)-1}{\CR(1,2,3,4)-1}.
\end{align*}
As $M_{0,n}$ is dense in $\Mbar_{0,n},$ these functional equations hold on $\Mbar_{0,n}$ wherever both sides take values in $\A^1.$

One may obtain global coordinates on  $M_{0,n}$ by composing forgetful morphisms with the cross-ratio map. Given $(C,p_1,\ldots,p_n)\in M_{0,n}$, we identify $C$ with $\P^1$ so that $p_1=\infty,$ $p_2=0$, and $p_3=1$; the induced map 
\begin{align*}
(\CR(1,2,3,4),\ldots,\CR(1,2,3,n)):M_{0,n}&\to \A^{n-3}\\
    (\P^1,\infty,0,1,p_4,\ldots,p_n)&\mapsto(p_4,\ldots,p_n)
\end{align*}
is an isomorphism onto its image $M_{0,n}\cong\A^{n-3}\setminus\{y_i=0,y_i=1,y_i=y_j\}.$ Via the cross-ratio functional equations, any cross-ratio $\CR(i_1,i_2,i_3,i_4)$ may be written as a rational function in the cross-ratios $\CR(1,2,3,i)$.

\subsection{Global coordinates on \texorpdfstring{$\H_{d,n}$}{H\_\{d,n\}}}\label{sec:HdnCoords}
We now work over $\Q(\zeta);$ in light of Remarks \ref{rem:BaseChange} and \ref{rem:BaseChange2}, we may now work with the enlarged marking set $\bA_{n}$ instead of $\bA_{n,0}$. For $a_{i_1},\ldots,a_{i_4}\in\mathbf{A}_n$ and $b_{i_1},\ldots,a_{i_4}\in\mathbf{B}_n$, there are cross-ratio maps $$\CR(C,a_{i_1},\ldots,a_{i_4}),\CR(D,b_{i_1},\ldots,b_{i_4}):\Hbar_{d,n}\P^1,$$ defined by pulling back from $M_{0,\bA_n}$ and $M_{0,\bB_n}$, respectively.

Given $(C,D,f)\in\H_{d,n}$, we choose coordinates on $C$ and $D$ such that $a_1=0,$ $a_{2,0}=1,$ $a_*=\infty$, $b_2=0,$ $b_3=1,$ $b_*=\infty$. With respect to these coordinates, $f:C\to D$ is identified with the function $z\mapsto z^d$, and we have $b_{i+1\mod n}=a_{i,k}^d$ for $i=2,\ldots,n$ and $k=0,\ldots,d-1$. In addition to the cross-ratio functional equations, we thus also have relations: 
\begin{align}\label{eq:RelsBtwnCR}
    \nonumber\CR(C,a_*,a_1,a_{2,0},a_{i,k})^d&=\CR(D,b_*,b_2,b_3,b_{(i+1\mod n)})\quad &i&=3,\ldots,n\\
    \CR(C,a_*,a_1,a_{2,0},a_{i,k})&=\zeta^k\CR(C,a_*,a_1,a_{2,0},a_{i,0})\quad &i&=3,\ldots,n, \quad k=1,\ldots,d-1\\\nonumber
    \CR(C,a_*,a_1,a_{2,0},a_{2,k})&=\zeta^k.
\end{align}
The tuple $(a_{3,0},\ldots, a_{n,0})=(\CR(C,a_*,a_1,a_{2,0},a_{3,0}),\ldots,\CR(C,a_*,a_1,a_{2,0},a_{n,0}))\in \A^{n-2}$ lies in the complement of the hypersurfaces defined by the equations $x_i=0,$ $x_i^d=1$, and $x_i^d=x_j^d$. Conversely, given a point in the complement of these hypersurfaces, one can reconstruct a point of $\H_{d,n}$. We conclude that $\H_{d,n}\cong\A^{n-2}\setminus\{x_i=0,x_i^d=1,x_i^d=x_j^d\}$. We can write down $\pi_1$ and $\pi_2$ with respect to these coordinates:
\begin{align*}
    \pi_1:\H_{d,n}\cong\A^{n-2}_{(x_i)}\setminus\{x_i=0,x_i^d=1,x_i^d=x_j^d\}&\to M_{0,n}\cong\A^{n-3}_{(y_i)}\setminus\{y_i=0,y_i=1,y_i=y_j\}\\
    (x_1,\ldots,x_{n-2})&\mapsto(\CR(\P^1,x_{n-2}^d,0,1,x_1^d),\ldots,\CR(\P^1,x_{n-2}^d,0,1,x_{n-3}^d))\\&\\
    \pi_2:\H_{d,n}\cong\A^{n-2}_{(x_i)}\setminus\{x_i=0,x_i^d=1,x_i^d=x_j^d\}&\to M_{0,n}\cong\A^{n-3}_{(z_i)}\setminus\{z_i=0,z_i=1,z_i=z_j\}\\
    (x_1,\ldots,x_{n-2})&\mapsto(\CR(\P^1,0,1,x_3,x_4),\ldots,\CR(\P^1,0,1,x_3,x_{n-2}))
\end{align*}

\subsection{Node smoothing parameters}\label{sec:NodeSmoothing}
  

Suppose $\sigma$ is an $n$-marked stable tree. For every edge $e$ of $\sigma$, we choose a $4$-tuple of legs as follows. If $v_{1},v_{2}\in\mathbf{V}(\sigma)$ are the endpoints of $e$, we choose distinct elements $i(e,1),\ldots,i(e,4)\in[n]$ satisfying
\begin{itemize}
\item $i(e,1)$ and $i(e,3)$ are in the connected component of $\sigma\setminus e$ that contains $v_{1}$ but does not contain $v_{2}$,
\item $i(e,1)$ and $i(e,3)$ are in different connected components of $\sigma\setminus \{v_{1}\}$,
\item $i(e,2)$ and $i(e,4)$ are in the connected component of $\sigma\setminus e$ that contains $v_{2}$ but does not contain $v_{1}$, and
\item $i(e,2)$ and $i(e,4)$ are in different connected components of $\sigma\setminus \{v_{2}\}$.
\end{itemize}
Such elements are guaranteed to exist by the stability condition of $\Mbar_{0,n}.$ Set $s_{e}:=\CR(i(e,1),\ldots,i(e,4))$. There is a Zariski neighborhood $W_\sigma$ of $Q_\sigma$ on which $s_e$ is regular for all $e\in\mathbf{E}(\sigma)$. Let $\sigma_e$ denote the stable tree obtained by contracting all edges in $\sigma$ except $e$. Then on $W_\sigma,$ we have $Z(s_e)=\bar{Q_{\sigma_e}}\cap W_{\sigma}$, and $Q_\sigma$ is the transverse intersection $\bigcap_{e\in\mathbf{E}(\sigma)}Z(s_e).$ The function $s_e$ is commonly referred to as a \textit{node smoothing parameter} for $e$. Indeed, the universal curve over $W_\sigma$ can be written, locally at the node corresponding to $e$, as the vanishing set of the polynomial $x_1x_2-s_e.$

More generally, we may define node smoothing parameters infinitesimally. If $\mathbf{y}\in Q_{\sigma}$, a node smoothing parameter for $e$ at $\mathbf{y}$ is an element $s\in\hat{\mathcal{O}}_{\mathbf{y},\Mbar_{0,n}}$ of the completion of the local ring of $\Mbar_{0,n}$ at $\mathbf{y}$ with the property that on $\Spec(\hat{\mathcal{O}}_{\mathbf{y},\Mbar_{0,n}})$, we have $Z(s)=\bar{Q_{\sigma_e}}$. A node smoothing parameter for $e$ (on $Q_\sigma$) is a function on the infinitesimal neighborhood of $Q_\sigma\subseteq\Mbar_{0,n}$ that restricts to a node smoothing parameter at every $\mathbf{y}\in Q_\sigma.$ Such functions exist for any $e\in\mathbf{E}(\sigma)$; an example is the cross-ratio function $s_e$, restricted to the infinitesimal neighborhood of $Q_\sigma$. We now record two observations that will be useful in Section \ref{sec:LocalAnalysis}.
\begin{obs}\label{obs:DifferentChoicesOfSmoothingParameter}
If $s_1$ and $s_2$ are two node smoothing parameters for $e$, then $s_1\sim s_2$ (in the notation from Section \ref{sec:Conventions}).
\end{obs}

\begin{obs}\label{obs:ProductOfSmoothingParameters}
Given $i_1,\ldots,i_4\in[n]$, let $\mathbf{E'}\subset\mathbf{E}$ be the (possibly empty) subset of edges $e$ with the property that $i_1$ and $i_3$ are in one connected component of $\sigma\setminus\{e\}$, and $i_2$ and $i_4$ are in the other. Then, locally at $Q_\sigma$, the locus $Z(\CR(i_1,\ldots,i_4))$ is the union of the hypersurfaces in $\{\bar{Q}_{\sigma_{e}}\}_{e\in\mathbf{E'}}$. In particular, if, for $e\in\mathbf{E'}$, $s_{e}$ is any choice of node smoothing parameter for $e$, then $\CR(i_1,\ldots,i_4)\sim \prod_{e\in \mathbf{E'}} s_{e}$.
\end{obs}

\subsection{Local coordinates on \texorpdfstring{$\Mbar_{0,n}$}{M\_\{0,n\}-bar}}\label{sec:M0nBarCoords}
Let $\sigma$ be an $n$-marked stable tree. For each vertex $v$ of $\sigma$ having valence four or more, we choose legs $i(v,1),\ldots, i(v,\val{v})\in[n]$ that are in $\val(v)$ distinct connected components of $\sigma\setminus\{v\}$. For $i=1,\ldots, \val(v)-3$, let $s_{v,i}=\CR(i(v,1),i(v,2),i(v,3),i(v,i+3))$. 
These parameters are called \textit{cross-ratio parameters} for the vertex $v$. Recall the product decomposition \eqref{eq:QSigmaProduct} of $Q_\sigma$; here $s_{v,i}$ is pulled back from $M_{0,\val(v)}$. It thus follows from Section \ref{sec:M0nCoords} that the parameters $((s_{v,i})_i)_{v}$ are global coordinates on $Q_\sigma,$ defining an open immersion $$Q_\sigma\cong\prod_{v}M_{0,\val(v)}\into\prod_{v\in\mathbf{V}(\sigma)}\A^{\val(v)-3}\cong\A^{n-3-\abs{\mathbf{E}(\sigma)}}.$$

Let $(s_e)_e$ be a choice of node smoothing parameters as in Section \ref{sec:NodeSmoothing}. Then by \cite{DeligneMumford1969,Knudsen1983}, the map $(s_{e})_{e}\times ((s_{v,i})_i)_{v}$ defines ``local coordinates on $\Mbar_{0,n}$ along $Q_\sigma$,'' in the sense that $(s_{e})_{e}\times ((s_{v,i})_i)_{v}$ is an isomorphism from the infinitesimal neighborhood of $Q_\sigma\subseteq\Mbar_{0,n}$ to the infinitesimal neighborhood of  $((s_{v,i})_i)_{v}(Q_\sigma)\subseteq\A^{n-3-\abs{\mathbf{E}(\sigma)}}\times\A^{\abs{\mathbf{E}(\sigma)}}\cong\A^{n-3}$. Algebraically, this infinitesimal neighborhood is $\Spec(A),$ where $$A=k[\{s_{v,i}^{\pm1},(s_{v,i}-1)^{-1},(s_{v,i}-s_{v,j})^{-1}\}_{v,i}][[\{s_e\}_e]].$$
If the node smoothing parameters $(s_e)_e$ are restrictions of local regular functions on $\Mbar_{0,n}$, this may be thought of more concretely; in this case, $(s_{e})_{e}\times ((s_{v,i})_i)_{v}$ is an \'etale map from a Zariski open neighborhood of $Q_\sigma$ to $\A^{n-3}$, that restricts to an isomorphism on $Q_\sigma$. Note, however, that we will need to do computations in the power series ring $A$; see Section \ref{sec:R6}.


\subsection{Local coordinates on \texorpdfstring{$\Hbar_{d,n}$}{H\_\{d,n\}-bar}}\label{sec:HnBarCoords}
Any 4-tuple of distinct elements $(a_{i_1},\ldots,a_{i_4})\subseteq \bA_n$ determines a cross-ratio map $\CR(C,a_{i_1},\ldots,a_{i_4}):\Hbar_{d,n}\to\P^1$; this map factors through the natural forgetful map $\Hbar_{d,n}\to\Mbar_{0,\bA_n}$. Similarly, any 4-tuple of distinct elements $(b_{i_1},\ldots,b_{i_4})\subseteq\bB_n$ determines a cross-ratio map $\CR(D,b_{i_1},\ldots,b_{i_4}):\Hbar_{d,n}\to\P^1$; this map factors through the natural forgetful map $\Hbar_{d,n}\to\Mbar_{0,\bB_n}$.
Since $\H_{d,n}$ is dense in $\Hbar_{d,n}$, each of the relations (\ref{eq:RelsBtwnCR}) holds wherever the relevant $\CR$ maps take values in $\A^1$, as do the cross-ratio functional equations pulled back from  $M_{0,\bA_n}$ and $M_{0,\bB_n}$.

Let $\gamma=(\sigma,\tau,\phi,\deg)$ be a combinatorial type as in Section \ref{sec:HnBar}, and let $R_\gamma\subseteq\Hbar_{d,n}$ be the corresponding locally closed boundary stratum. We now describe coordinates on the infinitesimal neighborhood $U_{\gamma}$ of $R_\gamma\subseteq\Hbar_{d,n}$, mimicking Section \ref{sec:M0nBarCoords}. 

For each vertex $w\in\mathbf{V}(\tau)$ with $\val(w)\ge4$, choose one $v\in\phi^{-1}(w)$. We choose an ordered $\val(w)$-tuple $(a_{v,1},\ldots,a_{v,\val(w)})$ of distinct elements of $\mathbf{A}_n$ with the property that $\phi(a_{v,1}),\ldots,\phi(a_{v,\val(w)})\in\mathbf{B}_n$ are in (the $\val(w)$) distinct connected components of $\tau\setminus\{w\}.$ For $i=1,\ldots, \val(w)-3$, consider the cross-ratio parameters $\tilde{s}_{w,i}:=\CR(C,a_{v,1},a_{v,2},a_{v,3},a_{v,i+3})$ for the vertex $v$. If $\phi^{-1}(w)=\{v\},$ then for convenience we choose the numbering so that the flags $\fl_1$ and $\fl_2$ connecting $v$ to $a_{v,1}$ and $a_{v,2}$ satisfy $\deg(\fl_1)=\deg(\fl_2)=d$. The product decomposition \eqref{eq:RGammaProduct} of $R_\gamma$ and Section \ref{sec:HdnCoords} imply that $((\tilde{s}_{w,i})_i)_w$ are global coordinates on $R_{\gamma}$, defining an open immersion
$$R_{\gamma}\cong \prod_{\substack{w\in\bV(\tau)\\\abs{\phi^{-1}(w)}=1}}\H_{d,\val(w)-1}\times \prod_{\substack{w\in\bV(\tau)\\ \abs{\phi^{-1}(w)}=d}}M_{0,\val(w)}\into\prod_{\substack{w\in\bV(\tau)\\\abs{\phi^{-1}(w)}=1}}\A^{\val(w)-3}\times \prod_{\substack{w\in\bV(\tau)\\ \abs{\phi^{-1}(w)}=d}}\A^{\val(w)-3}\cong\A^{n-2-\abs{\mathbf{E}(\tau)}}.$$

For each edge $\eta\in\mathbf{E}(\tau),$ we choose an edge $e\in\phi^{-1}(\eta)$ and a node smoothing parameter for $e$ as in Section \ref{sec:NodeSmoothing}. It follows from \cite[page 61]{HarrisMumford1982} that $(\tilde{s}_{\eta})_{\eta}\times((\tilde{s}_{w,i})_i)_w$ defines an isomorphism from the infinitesimal neighborhood $U_\gamma$ of $R_\gamma\subseteq\Hbar_{d,n}$ to the infinitesimal neighborhood of $((\tilde{s}_{w,i})_i)_w(R_\gamma)\subseteq\A^{n-2-\abs{\mathbf{E}(\tau)}}\times\A^{\abs{\mathbf{E}(\tau)}}\cong\A^{n-2}$. Algebraically, we have $U_\gamma\cong\Spec(B),$ where $$B=k[\{\tilde{s}_{w,i}^{\pm1},(\tilde{s}_{w,i}\pm1)^{-1},(\tilde{s}_{w,i}\pm\tilde{s}_{w,j})^{-1}\}_{\abs{\phi^{-1}(w)}=1},\{\tilde{s}_{w,i}^{\pm1},(\tilde{s}_{w,i}-1)^{-1},(\tilde{s}_{w,i}-\tilde{s}_{w,j})^{-1}\}_{\abs{\phi^{-1}(w)}=d}][[\{\tilde{s}_\eta\}_\eta]].$$ In $U_\gamma$, $R_\gamma$ is the vanishing set of the coordinates $(\tilde{s}_{\eta})_{\eta}$. As in Section \ref{sec:M0nBarCoords}, if the node smoothing parameters $\tilde{s}_\eta$ are restrictions of local regular functions on $\Hbar_{d,n}$, then $(\tilde{s}_{\eta})_{\eta}\times((\tilde{s}_{w,i})_i)_w$ is an \'etale map from a Zariski open neighborhood of $R_\gamma$ to $\A^{n-2}$, that restricts to an isomorphism on $R_\gamma.$

\subsection{Local equations for \texorpdfstring{$\pi_1$}{pi\_1}, \texorpdfstring{$\pi_2$}{pi\_2} and \texorpdfstring{$\tilde{\Delta}_{d,n}$}{Delta\_\{d,n\}-tilde}}\label{sec:pisanddeltacoords}

Let $\gamma$, $R_\gamma$ $U_\gamma$, and $(\tilde{s}_{\eta})_{\eta}\times((\tilde{s}_{w,i})_i)_w$ be as in Section \ref{sec:HnBarCoords}. 
Any regular function on $U_\gamma$ may be expressed in terms of the coordinates $(\tilde{s}_{\eta})_{\eta}\times((\tilde{s}_{w,i})_i)_w$. 

More concretely, assume the node smoothing parameters $\tilde{s}_{\eta}$ are cross-ratio functions pulled back from $M_{0,\bA_n}$ (as in Section \ref{sec:NodeSmoothing}). Then any function of the form $\CR(C,a_{i_1},\ldots,a_{i_4})$ for $a_{i_j}\in\bA_n$ or $\CR(D,b_{i_1},\ldots,b_{i_4})$ for $b_{i_j}\in\bB_n$ can be expressed as rational functions in the chosen coordinates using the cross-ratio functional equations as well as the relations (\ref{eq:RelsBtwnCR}), as follows. A cross-ratio $F:=\CR(D,b_{i_1},\ldots,b_{i_4})$ can be written, using the cross-ratio functional equations pulled back from $M_{0,\bB_n}$, as a rational function $F'$ in the cross-ratios $\CR(b_*,b_2,b_3,b_i)$. The first equation in \eqref{eq:RelsBtwnCR} then expresses $F'$ as a polynomial $F''$ in cross-ratios pulled back from $M_{0,\bA_n}$. Finally, the cross-ratio functional equations pulled back from $M_{0,\bA_n}$ are used to express $F''$ as an expression $F'''$ in terms of the chosen coordinates.
\begin{caution}\label{cau:Roots}
The map $(\tilde{s}_{\eta})_{\eta}\times((\tilde{s}_{w,i})_i)_w$ is \'etale in a Zariski neighborhood of $R_\gamma$, but possibly finite-to-one. Thus the last sentence of the previous paragraph may require solving polynomial equations in the coordinates --- so $F'''$ is not necessarily a rational function in the coordinates, but an expression involving roots. (See Section \ref{sec:R6}.) Since $(\tilde{s}_{\eta})_{\eta}\times((\tilde{s}_{w,i})_i)_w$ restricts to an isomorphism on $U_\gamma$, $F'''$ is a well-defined element of the ring $B$ above. However, note that $(\tilde{s}_{\eta})_{\eta}\times((\tilde{s}_{w,i})_i)_w$ may map other boundary strata to the image of $R_\gamma$; one must make sure to take the expansion of $F'''$ along $R\gamma,$ and not another such stratum. Practically speaking, $F'''$ is an expression in terms of roots, and one must carefully choose the correct branch of each root to obtain the correct function on $U_\gamma$.
\end{caution}

The process above allows us to write $\pi_1|_{U_\gamma}$ and $\pi_2|_{U_{\gamma}}$ in coordinates. Let $\bar{\tau}$ and $\bar{\sigma}$ be as in Section \ref{sec:HnBar}. The image $\pi_1(R_{\gamma})$ is the locally closed stratum $Q_{\bar{\tau}}\subseteq\Mbar_{0,n}$, and the image $\pi_2(R_{\gamma})$ is contained in the locally closed stratum $Q_{\bar{\sigma}}\subseteq\Mbar_{0,n}$. Then $\pi_2^*(x_i)$ is of the form $\CR(C,a_{i_1},\ldots, a_{i_4})$, and $\pi_1^{*}(y_i)$ is of the form $\CR(D,b_{i_1},\ldots, b_{i_4})$. Thus, using cross-ratio functional equations and the relations (\ref{eq:RelsBtwnCR}), each $\pi_1^{*}(y_i)$ and $\pi_2^{*}(x_i)$ can be written as rational functions in the parameters $(\tilde{s}_{\eta})_{\eta}\times((\tilde{s}_{w,i})_i)_w$.

Suppose $\gamma=(\sigma,\tau,\phi,\deg)$ is such that $\bar{\sigma}=\bar{\tau}$. (By Observation \ref{obs:CombinatorialCondition}, this is true whenever $R_\gamma\cap\tilde{\Delta}_{d,n}\ne\emptyset.$) We can pick local coordinates $x_1,\ldots,x_{n-3}$ on the infinitesimal neighborhood of $Q_{\bar{\sigma}}=Q_{\bar{\tau}}\subseteq\Mbar_{0,n}$.
Then on $U_\gamma$, $\tilde{\Delta}_{d,n}$ is cut out by the $n-3$ equations $\pi_1^{*}(x_1)=\pi_2^{*}(x_1),\ldots,\pi_1^{*}(x_{n-3})=\pi_2^{*}(x_{n-3})$. As above, each of these equations may be written in terms of the coordinates $(\tilde{s}_{\eta})_{\eta}\times((\tilde{s}_{w,i})_i)_w$. Thus the problem of determining $R_\gamma\cap\tilde{\Delta}_{d,n}$ is reduced to an algebraic computation in the ring $B.$ (See Sections \ref{sec:LocalAnalysis} and \ref{sec:LocalAnalysis2} for computations demonstrating this process.) In summary, we have:
\begin{algorithm}\label{alg:Punctures}
For a fixed $d$ and $n$,
\begin{itemize}
    \item We enumerate combinatorial types $\gamma$ as in Section \ref{sec:HnBar}.
    \item For each such $\gamma,$ we choose coordinates for an \'etale neighborhood $U_\gamma$ of $R_\gamma$ as in Section \ref{sec:HnBarCoords}.
    \item We give defining equations for $\tilde{\Delta}_{d,n}$ in $U_\gamma$ as in Section \ref{sec:pisanddeltacoords}.
    \item We compute the irreducible component $\overline{\Per_{d,n}}\subseteq\tilde{\Delta}_{d,n}$ in $U_\gamma$ and the intersection $\overline{\Per_{d,n}}\cap R_\gamma$ via algebraic computations (e.g. the former by computing a primary decomposition).
\end{itemize}
Altogether, this gives local equations for $\Per_{d,n}$ at all punctures.
\end{algorithm}
\begin{rem}
Observations \ref{obs:CombinatorialCondition}, \ref{obs:DifferentChoicesOfSmoothingParameter}, \ref{obs:ProductOfSmoothingParameters}, as well as Observations \ref{obs:SmoothingParametersOnSourceAndTarget} and \ref{obs:CrossRatiosParametersOnSourceAndTarget} below, provide significant shortcuts to many of the calculations.
\end{rem}
\begin{rem}\label{rem:Generalize}
Algorithm \ref{alg:Punctures} can be adapted to arbitrary COR varieties, with is one significant change which we now describe. In Section \ref{sec:HnBarCoords}, the cross-ratio and node smoothing parameters are algebraically independent, and define local coordinates near the stratum $R_\gamma$. In general, however, there may be relations among these parameters, see Section 4 of \cite{HarrisMumford1982}.
\end{rem}

Note that no Thurston/Epstein-type transversality result holds on the boundary of $\Hbar_{d,n}$; $\tilde{\Delta}_{d,n}$ need not be smooth and may have components of dimension greater than 1 (the expected dimension). (In Section \ref{sec:C5InHBar}, we see that $\tilde{\Delta}_{2,5}$ is not smooth.)

We record two more observations that will simplify our calculations in Section \ref{sec:LocalAnalysis}.
\begin{obs}\label{obs:SmoothingParametersOnSourceAndTarget}
 Let $\eta\in\mathbf{E}(\tau)$ and let $e\in\phi^{-1}(\eta).$ Let $\tilde{s},s':\Hbar_{d,n}\to\P^1$ be node smoothing parameters for $e$ and $\eta$, respectively. On the infinitesimal neighborhood $U_\gamma$ of $R_\gamma\subseteq\Hbar_{d,n},$ $\tilde{s}$ and $s'$ are regular functions, and the balancing condition at nodes (Section \ref{sec:HnBar}) implies that $\tilde{s}\sim\pi_2^*((s')^{\deg(e)})$.
\end{obs}
\begin{obs}\label{obs:CrossRatiosParametersOnSourceAndTarget}
Let $w\in\mathbf{V}(\tau)$, and let $v\in\phi^{-1}(w)$. Suppose $\deg(v)=1.$ Let $s_w$ be a cross-ratio parameter for $w$ defined by a 4-tuple of flags $(\fl_1,\ldots,\fl_4)$ at $w.$ Let $\tilde{s}_v$ be a cross-ratio parameter for $v$ defined by a 4-tuple of flags $(\tilde{\fl}_1,\ldots,\tilde{\fl}_4)$ at $v,$ and suppose $\phi(\tilde{\fl}_i)=\fl_i$ for $i\in\{1,\ldots,4\}.$ Then on $R_\gamma$, $\pi_2^*(s_w)=\pi_1^*(\tilde{s}_v).$ Equivalently, as functions on $U_\gamma$, the function $\pi_2^*(s_w)-\pi_1^*(\tilde{s}_v)$ is in the ideal generated by the node smoothing parameters $(\tilde{s}_\eta)_{\eta}$.

Suppose instead $\deg(v)=d$, and further that $\deg(\fl_1)=\deg(\fl_2)=d$. Then choosing coordinates on the irreducible component corresponding to $v$, we see that $\pi_2^*(s_w)=\pi_1^*(\tilde{s}_v)^d$ on $R_\gamma$.
\end{obs}

\section{Application 1: The geometry of \texorpdfstring{$\overline{\Per_{2,5}}$}{Per\_\{2,5\}} as a curve in \texorpdfstring{$\Hbar_{2,5}$}{H\_\{2,5\}-bar}}\label{sec:C5InHBar}

We now restrict our attention to the case $(d,n)=(2,5)$. Using the boundary stratification of $\Hbar_{2,5},$ we will analyze the local geometry of $\tilde{\Delta}_{2,5}\setminus\Per_{2,5}^{\circ}$. This will allow us to identify the points of $\overline{\Per_{2,5}}\setminus\Per_{2,5}^{\circ}$. The analysis will also show that all of these boundary points are smooth points of $\overline{\Per_{2,5}}$; since $\Per_{2,5}^{\circ}$ is smooth, this will imply that $\overline{\Per_{2,5}}$ is smooth and $\widetilde{\Per}_{2,5}=\overline{\Per_{2,5}}.$

\subsection{Enumeration of boundary strata}\label{sec:StrataEnumeration}
In order to identify $\overline{\Per_{2,5}}\setminus\Per_{2,5}^{\circ}$, we first describe $\tilde{\Delta}_{2,5}\setminus \H_{2,5}$. We now list all boundary strata $R_\gamma$ in $\Hbar_{2,5}$ that intersect $\tilde{\Delta}_{2,5}$. Observation \ref{obs:CombinatorialCondition} gives a combinatorial necessary condition on $\gamma$; we have enumerated all $\gamma$ that satisfy this condition, and checked our answer by implementing the condition on a computer. The results appear in Tables \ref{fig:StrataThatIntersectC5} and \ref{fig:DoNotIntersectC5}, according to whether or not $\gamma$ is found to intersect $\overline{\Per_{2,5}}$ in Section \ref{sec:LocalAnalysis}.

\begin{table}\centering
  
  \begin{tabular}{ccc}
    \begin{tikzpicture}
      \draw (-1.5,0) node {$\tau$};
\draw (-1.5,2) node {$\sigma$};
\draw
      (0,0)--(3,0);
\draw (0,0) node {$\bullet$};
\draw
      (0,0)--++(160:.4);
\draw (0,0)++(160:.4) node[left] {2};
\draw
      (0,0)--++(200:.4);
\draw (0,0)++(200:.4) node[left] {$*$};
\draw
      (1,0) node {$\bullet$};
\draw (1,0)--++(-90:.4);
\draw
      (1,0)++(-90:.4) node[below] {1};
\draw (2,0) node {$\bullet$};
      \draw (2,0)--++(-90:.4);
\draw (2,0)++(-90:.4) node[below] {5};
      \draw (3,0) node {$\bullet$};
\draw (3,0)--++(20:.4);
\draw
      (3,0)++(20:.4) node[right] {3};
\draw (3,0)--++(-20:.4);
\draw
      (3,0)++(-20:.4) node[right] {4};
\draw[->]
      (1.5,.75)--(1.5,0.25);
\draw
      (3,2.5)--(2,2.5)--(1,2.5)--(0,2)--(1,1.5)--(2,1.5)--(3,1.5);
      \draw (3,2.5) node {$\bullet$};
\draw (3,2.5)--++(20:.4);
\draw
      (3,2.5)++(20:.4) node[right] {2};
\draw (2,2.5) node
      {$\bullet$};
\draw (1,2.5) node {$\bullet$};
\draw (0,2) node
      {$\bullet$};
\draw[very thick] (0,2)--++(160:.4);
\draw
      (0,2)++(160:.4) node[left] {1};
\draw[very thick]
      (0,2)--++(200:.4);
\draw (0,2)++(200:.4) node[left] {$*$};
\draw (0.5,1.75) node[below] {\footnotesize$s_1$};
\draw
      (1,1.5) node {$\bullet$};
\draw (1,1.5)--++(-90:.4);
\draw
      (1,1.5)++(-90:.4) node[below] {5};
      \draw (1.5,1.5) node[above] {\footnotesize$s_2$};
\draw (2,1.5) node
      {$\bullet$};
\draw (2,1.5)--++(-90:.4);
\draw (2,1.5)++(-90:.4)
      node[below] {4};
      \draw (2.5,1.5) node[above] {\footnotesize$s_3$};
\draw (3,1.5) node {$\bullet$};
\draw
      (3,1.5)--++(-20:.4);
\draw (3,1.5)++(-20:.4) node[right] {3};
    \end{tikzpicture}
    &
      \begin{tikzpicture}
        \draw (0,0)--(3,0);
\draw (0,0) node {$\bullet$};
\draw
        (0,0)--++(160:.4);
\draw (0,0)++(160:.4) node[left] {2};
\draw
        (0,0)--++(200:.4);
\draw (0,0)++(200:.4) node[left] {$*$};
        \draw (1,0) node {$\bullet$};
\draw (1,0)--++(-90:.4);
\draw
        (1,0)++(-90:.4) node[below] {3};
\draw (2,0) node {$\bullet$};
        \draw (2,0)--++(-90:.4);
\draw (2,0)++(-90:.4) node[below]
        {1};
\draw (3,0) node {$\bullet$};
\draw (3,0)--++(-20:.4);
        \draw (3,0)++(-20:.4) node[right] {4};
\draw (3,0)--++(20:.4);
        \draw (3,0)++(20:.4) node[right] {5};
\draw[->]
        (1.5,.75)--(1.5,0.25);
\draw
        (3,2.5)--(2,2.5)--(1,2.5)--(0,2)--(1,1.5)--(2,1.5)--(3,1.5);
        \draw (3,2.5) node {$\bullet$};
\draw (3,2.5)--++(20:.4);
        \draw (3,2.5)++(20:.4) node[right] {4};
\draw (2,2.5) node
        {$\bullet$};
\draw (2,2.5)--++(-90:.4);
\draw
        (2,2.5)++(-90:.4) node[below] {5};
\draw (1,2.5) node
        {$\bullet$};
\draw (0,2) node {$\bullet$};
\draw[very thick]
        (0,2)--++(160:.4);
\draw (0,2)++(160:.4) node[left] {1};
        \draw[very thick] (0,2)--++(200:.4);
\draw (0,2)++(200:.4)
        node[left] {$*$};
        \draw (0.5,1.75) node[below] {\footnotesize$s_1$};
\draw (1,1.5) node {$\bullet$};
\draw
        (1,1.5)--++(-90:.4);
\draw (1,1.5)++(-90:.4) node[below] {2};
\draw (1.5,1.5) node[below] {\footnotesize$s_2$};
        \draw (2,1.5) node {$\bullet$};
        \draw (2.5,1.5) node[below] {\footnotesize$s_3$};
\draw (3,1.5) node
        {$\bullet$};
\draw (3,1.5)--++(-20:.4);
\draw
        (3,1.5)++(-20:.4) node[right] {3};
      \end{tikzpicture}
    &
      \begin{tikzpicture}
        \draw (0,0)--(2,0);
\draw (1,0)--++(-45:.6);
\draw (0,0) node
        {$\bullet$};
\draw (0,0)--++(160:.4);
\draw (0,0)++(160:.4)
        node[left] {2};
\draw (0,0)--++(200:.4);
\draw (0,0)++(200:.4)
        node[left] {$*$};
\draw (1,0) node {$\bullet$};
\draw
        (1,0)++(-45:.6) node {$\bullet$};
\draw
        (1,0)++(-45:.6)--++(-110:.4);
\draw (1,0)++(-45:.6)++(-110:.4)
        node[below] {1};
\draw (1,0)++(-45:.6)--++(-70:.4);
\draw
        (1,0)++(-45:.6)++(-70:.4) node[below] {4};
\draw (2,0) node
        {$\bullet$};
\draw (2,0)--++(20:.4);
\draw (2,0)++(20:.4)
        node[right] {5};
\draw (2,0)--++(-20:.4);
\draw
        (2,0)++(-20:.4) node[right] {3};
\draw[->] (1,.75)--(1,0.25);
        \draw (2,2.5)--(1,2.5)--(0,2)--(1,1.5)--(2,1.5);
\draw
        (1,2.5)--++(-45:.6);
\draw (1,1.5)--++(-45:.6);
\draw (2,2.5)
        node {$\bullet$};
\draw (2,2.5)--++(20:.4);
\draw
        (2,2.5)++(20:.4) node[right] {4};
\draw (1,2.5) node
        {$\bullet$};
\draw (1,2.5)++(-45:.6) node {$\bullet$};
\draw
        (0,2) node {$\bullet$};
\draw[very thick] (0,2)--++(160:.4);
        \draw (0,2)++(160:.4) node[left] {1};
\draw[very thick]
        (0,2)--++(200:.4);
\draw (0,2)++(200:.4) node[left] {$*$};
        \draw (1,1.5) node {$\bullet$};
\draw (1,1.5)++(-45:.6) node
        {$\bullet$};
\draw (1,1.5)++(-45:.6)--++(-110:.4);
\draw
        (1,1.5)++(-45:.6)++(-110:.4) node[below] {5};
\draw
        (1,1.5)++(-45:.6)--++(-70:.4);
\draw
        (1,1.5)++(-45:.6)++(-70:.4) node[below] {3};
\draw (2,1.5)
        node {$\bullet$};
\draw (2,1.5)--++(-20:.4);
\draw
        (2,1.5)++(-20:.4) node[right] {2};
      \end{tikzpicture}
    \\$(\gamma_1)$&$(\gamma_2)$&$(\gamma_3)$\\
    &&\\\hline
    \begin{tikzpicture}
      \draw (-1.5,0) node {$\tau$};
\draw (-1.5,2) node
      {$\sigma$};
\draw (0,0)--(2,0);
\draw (0,0) node
      {$\bullet$};
\draw (0,0)--++(160:.4);
\draw (0,0)++(160:.4)
      node[left] {2};
\draw (0,0)--++(200:.4);
\draw (0,0)++(200:.4)
      node[left] {$*$};
\draw (1,0) node
      {$\bullet$};
\draw (1,0)--++(-110:.4);
\draw (1,0)++(-110:.4)
      node[below] {3};
\draw (1,0)--++(-70:.4);
\draw (1,0)++(-70:.4)
      node[below] {4};
\draw (2,0) node
      {$\bullet$};
\draw (2,0)--++(20:.4);
\draw (2,0)++(20:.4)
      node[right] {1};
\draw (2,0)--++(-20:.4);
\draw (2,0)++(-20:.4)
      node[right] {5};
\draw[->] (1,.75)--(1,0.25);
\draw
      (2,2.5)--(1,2.5)--(0,2)--(1,1.5)--(2,1.5);
\draw (2,2.5) node
      {$\bullet$};
\draw (2,2.5)--++(20:.4);
\draw (2,2.5)++(20:.4)
      node[right] {5};
\draw (1,2.5) node
      {$\bullet$};
\draw (0,2) node
      {$\bullet$};
\draw[very thick] (0,2)--++(160:.4);
\draw
      (0,2)++(160:.4) node[left] {1};
\draw[very thick]
      (0,2)--++(200:.4);
\draw (0,2)++(200:.4) node[left]
      {$*$};
      \draw (0.5,1.75) node[below] {\footnotesize$s_1$};
\draw (1,1.5) node
      {$\bullet$};
\draw (1,1.5)--++(-110:.4);
\draw
      (1,1.5)++(-110:.4) node[below] {2};
\draw (1,1.5)--++(-70:.4);
      \draw (1,1.5)++(-70:.4) node[below] {3};
      \draw (1.5,1.5) node[above] {\footnotesize$s_2$};
\draw (2,1.5) node
      {$\bullet$};
\draw (2,1.5)--++(-20:.4);
\draw (2,1.5)++(-20:.4)
      node[right] {4};
    \end{tikzpicture}
    &
      \begin{tikzpicture}
        \draw (0,0)--(1,0);
\draw (0,0) node
        {$\bullet$};
\draw (0,0)--++(140:.4);
\draw (0,0)++(140:.4)
        node[left] {2};
\draw (0,0)--++(220:.4);
\draw (0,0)++(220:.4)
        node[left] {4};
\draw (0,0)--++(180:.4);
\draw (0,0)++(180:.4)
        node[left] {$*$};
\draw (1,0) node
        {$\bullet$};
\draw (1,0)--++(-40:.4);
\draw (1,0)++(-40:.4)
        node[right] {5};
\draw (1,0)--++(0:.4);
\draw (1,0)++(0:.4)
        node[right] {3};
\draw (1,0)--++(40:.4);
\draw (1,0)++(40:.4)
        node[right] {1};
\draw[->] (.5,1)--(.5,0.5);
\draw
        (1,2.5)--(0,2)--(1,1.5);
\draw (1,2.5) node
        {$\bullet$};
\draw (1,2.5)--++(40:.4);
\draw (1,2.5)++(40:.4)
        node[right] {5};
\draw (0,2) node
        {$\bullet$};
\draw[very thick] (0,2)--++(140:.4);
\draw
        (0,2)++(140:.4) node[left] {1};
\draw (0,2)--++(220:.4);
\draw
        (0,2)++(220:.4) node[left] {3};
\draw[very thick]
        (0,2)--++(180:.4);
\draw (0,2)++(180:.4) node[left]
        {$*$};
        \draw (0.5,1.75) node[below] {\footnotesize$s_1$};
\draw (1,1.5) node
        {$\bullet$};
\draw (1,1.5)--++(-40:.4);
\draw
        (1,1.5)++(-40:.4) node[right] {4};
\draw (1,1.5)--++(0:.4);
        \draw (1,1.5)++(0:.4) node[right] {2};
      \end{tikzpicture}
    &
      \begin{tikzpicture}
        \draw (0,0)--(1,0);
\draw (0,0) node
        {$\bullet$};
\draw (0,0)--++(140:.4);
\draw (0,0)++(140:.4)
        node[left] {2};
\draw (0,0)--++(220:.4);
\draw (0,0)++(220:.4)
        node[left] {5};
\draw (0,0)--++(180:.4);
\draw (0,0)++(180:.4)
        node[left] {$*$};
\draw (1,0) node
        {$\bullet$};
\draw (1,0)--++(-40:.4);
\draw (1,0)++(-40:.4)
        node[right] {4};
\draw (1,0)--++(0:.4);
\draw (1,0)++(0:.4)
        node[right] {3};
\draw (1,0)--++(40:.4);
\draw (1,0)++(40:.4)
        node[right] {1};
\draw[->] (.5,1)--(.5,0.5);
\draw
        (1,2.5)--(0,2)--(1,1.5);
\draw (1,2.5) node
        {$\bullet$};
\draw (1,2.5)--++(40:.4);
\draw (1,2.5)++(40:.4)
        node[right] {5};
\draw (1,2.5)--++(0:.4);
\draw
        (1,2.5)++(0:.4) node[right] {2};
        \draw (0.5,2.25) node[above] {\footnotesize$s_1$};
\draw (0,2) node
        {$\bullet$};
\draw[very thick] (0,2)--++(140:.4);
\draw
        (0,2)++(140:.4) node[left] {1};
\draw (0,2)--++(220:.4);
\draw
        (0,2)++(220:.4) node[left] {4};
\draw[very thick]
        (0,2)--++(180:.4);
\draw (0,2)++(180:.4) node[left]
        {$*$};
\draw (1,1.5) node
        {$\bullet$};
\draw (1,1.5)--++(-40:.4);
\draw
        (1,1.5)++(-40:.4) node[right] {3};
      \end{tikzpicture}
    \\$(\gamma_4)$&$(\gamma_5)$&$(\gamma_6)$\\
    &&\\\hline
    \begin{tikzpicture}
      \draw (-1.5,0) node {$\tau$};
\draw (-1.5,2) node {$\sigma$};
\draw
      (0,0)--(1,0);
\draw (0,0) node {$\bullet$};
\draw
      (0,0)--++(160:.4);
\draw (0,0)++(160:.4) node[left] {2};
\draw
      (0,0)--++(200:.4);
\draw (0,0)++(200:.4) node[left] {$*$};
\draw
      (1,0) node {$\bullet$};
\draw (1,0)--++(-60:.4);
\draw
      (1,0)++(-60:.4) node[right] {3};
\draw (1,0)--++(-20:.4);
\draw
      (1,0)++(-20:.4) node[right] {4};
\draw (1,0)--++(20:.4);
\draw
      (1,0)++(20:.4) node[right] {5};
\draw (1,0)--++(60:.4);
\draw
      (1,0)++(60:.4) node[right] {1};
\draw[->] (.5,1)--(.5,0.5);
      \draw (1,2.5)--(0,2)--(1,1.5);
\draw (1,2.5) node {$\bullet$};
      \draw (0,2) node {$\bullet$};
\draw[very thick]
      (0,2)--++(160:.4);
\draw (0,2)++(160:.4) node[left] {1};
      \draw[very thick] (0,2)--++(200:.4);
\draw (0,2)++(200:.4)
      node[left] {$*$};
\draw (1,1.5) node {$\bullet$};
\draw
      (1,1.5)--++(-60:.4);
\draw (1,1.5)++(-60:.4) node[right] {2};
      \draw (1,1.5)--++(-20:.4);
\draw (1,1.5)++(-20:.4) node[right]
      {3};
\draw (1,1.5)--++(20:.4);
\draw (1,1.5)++(20:.4)
      node[right] {4};
\draw (1,1.5)--++(60:.4);
\draw
      (1,1.5)++(60:.4) node[right] {5};
    \end{tikzpicture}&
                       \begin{tikzpicture}
                         \draw (0,0)--(1,0);
\draw (0,0) node
                         {$\bullet$};
\draw (0,0)--++(160:.4);
\draw
                         (0,0)++(160:.4) node[left] {1};
\draw
                         (0,0)--++(200:.4);
\draw (0,0)++(200:.4)
                         node[left] {$*$};
\draw (1,0) node
                         {$\bullet$};
\draw (1,0)--++(-60:.4);
\draw
                         (1,0)++(-60:.4) node[right] {2};
\draw
                         (1,0)--++(-20:.4);
\draw (1,0)++(-20:.4)
                         node[right] {3};
\draw (1,0)--++(20:.4);
                         \draw (1,0)++(20:.4) node[right] {4};
\draw
                         (1,0)--++(60:.4);
\draw (1,0)++(60:.4)
                         node[right] {5};
\draw[->] (.5,1)--(.5,0.5);
                         \draw[very thick] (1,1.5)--(0,1.5);
\draw (0,1.5) node
                         {$\bullet$};
\draw (0,1.5)--++(160:.4);
\draw
                         (0,1.5)++(160:.4) node[left] {5};
\draw[very
                         thick] (0,1.5)--++(200:.4);
\draw
                         (0,1.5)++(200:.4) node[left] {$*$};
\draw
                         (1,1.5) node {$\bullet$};
\draw[very thick]
                         (1,1.5)--++(-60:.4);
\draw (1,1.5)++(-60:.4)
                         node[right] {1};
\draw (1,1.5)--++(-20:.4);
                         \draw (1,1.5)++(-20:.4) node[right] {2};
                         \draw (1,1.5)--++(20:.4);
\draw
                         (1,1.5)++(20:.4) node[right] {3};
\draw
                         (1,1.5)--++(60:.4);
\draw (1,1.5)++(60:.4)
                         node[right] {4};
                       \end{tikzpicture}
    &

     \begin{tikzpicture}
       \draw (0,0)--(1,0);
\draw (0,0) node {$\bullet$};
\draw
       (0,0)--++(160:.4);
\draw (0,0)++(160:.4) node[left] {3};
\draw
       (0,0)--++(200:.4);
\draw (0,0)++(200:.4) node[left] {$*$};
       \draw (1,0) node {$\bullet$};
\draw (1,0)--++(-60:.4);
\draw
       (1,0)++(-60:.4) node[right] {1};
\draw (1,0)--++(-20:.4);
\draw
       (1,0)++(-20:.4) node[right] {2};
\draw (1,0)--++(20:.4);
\draw
       (1,0)++(20:.4) node[right] {4};
\draw (1,0)--++(60:.4);
\draw
       (1,0)++(60:.4) node[right] {5};
\draw[->] (.5,1)--(.5,0.5);
       \draw[very thick] (1,1.5)--(0,1.5);
\draw (0,1.5) node {$\bullet$};
\draw
       (0,1.5)--++(160:.4);
\draw (0,1.5)++(160:.4) node[left] {2};
       \draw[very thick] (0,1.5)--++(200:.4);
\draw (0,1.5)++(200:.4)
       node[left] {$*$};
\draw (1,1.5) node {$\bullet$};
\draw
       (1,1.5)--++(-60:.4);
\draw (1,1.5)++(-60:.4) node[right] {5};
       \draw[very thick] (1,1.5)--++(-20:.4);
\draw (1,1.5)++(-20:.4)
       node[right] {1};
\draw (1,1.5)--++(20:.4);
\draw
       (1,1.5)++(20:.4) node[right] {3};
\draw (1,1.5)--++(60:.4);
       \draw (1,1.5)++(60:.4) node[right] {4};
     \end{tikzpicture}
    \\$(\gamma_7)$&$(\gamma_{\mathrm{I}})$&$(\gamma_{\mathrm{II}})$\\
    &&\\\hline
    \begin{tikzpicture}
      \draw (-1.5,0) node {$\tau$};
\draw (-1.5,1.5) node
      {$\sigma$};
\draw (0,0)--(1,0);
\draw (0,0) node
      {$\bullet$};
\draw (0,0)--++(160:.4);
\draw (0,0)++(160:.4)
      node[left] {4};
\draw (0,0)--++(200:.4);
\draw (0,0)++(200:.4)
      node[left] {$*$};
\draw (1,0) node
      {$\bullet$};
\draw (1,0)--++(-60:.4);
\draw (1,0)++(-60:.4)
      node[right] {1};
\draw (1,0)--++(-20:.4);
\draw (1,0)++(-20:.4)
      node[right] {2};
\draw (1,0)--++(20:.4);
\draw (1,0)++(20:.4)
      node[right] {3};
\draw (1,0)--++(60:.4);
\draw (1,0)++(60:.4)
      node[right] {5};
\draw[->] (.5,1)--(.5,0.5);
\draw[very thick]
      (1,1.5)--(0,1.5);
\draw (0,1.5) node
      {$\bullet$};
\draw (0,1.5)--++(160:.4);
\draw (0,1.5)++(160:.4)
      node[left] {3};
\draw[very thick] (0,1.5)--++(200:.4);
\draw
      (0,1.5)++(200:.4) node[left]
      {$*$};
\draw (1,1.5) node
      {$\bullet$};
\draw (1,1.5)--++(-60:.4);
\draw (1,1.5)++(-60:.4)
      node[right] {5};
\draw[very thick] (1,1.5)--++(-20:.4);
\draw
      (1,1.5)++(-20:.4) node[right] {1};
\draw (1,1.5)--++(20:.4);
      \draw (1,1.5)++(20:.4) node[right] {2};
\draw
      (1,1.5)--++(60:.4);
\draw (1,1.5)++(60:.4) node[right] {4};
    \end{tikzpicture}
    &
      \begin{tikzpicture}
        \draw (0,0)--(1,0);
\draw (0,0) node
        {$\bullet$};
\draw (0,0)--++(160:.4);
\draw (0,0)++(160:.4)
        node[left] {5};
\draw (0,0)--++(200:.4);
\draw (0,0)++(200:.4)
        node[left] {$*$};
\draw (1,0) node
        {$\bullet$};
\draw (1,0)--++(-60:.4);
\draw (1,0)++(-60:.4)
        node[right] {1};
\draw (1,0)--++(-20:.4);
\draw
        (1,0)++(-20:.4) node[right] {2};
\draw (1,0)--++(20:.4);
\draw
        (1,0)++(20:.4) node[right] {3};
\draw (1,0)--++(60:.4);
\draw
        (1,0)++(60:.4) node[right] {4};
\draw[->] (.5,1)--(.5,0.5);
        \draw[very thick] (1,1.5)--(0,1.5);
\draw (0,1.5) node
        {$\bullet$};
\draw (0,1.5)--++(160:.4);
\draw
        (0,1.5)++(160:.4) node[left] {4};
\draw[very thick]
        (0,1.5)--++(200:.4);
\draw (0,1.5)++(200:.4) node[left]
        {$*$};
\draw (1,1.5) node
        {$\bullet$};
\draw (1,1.5)--++(-60:.4);
\draw
        (1,1.5)++(-60:.4) node[right] {5};
\draw[very thick]
        (1,1.5)--++(-20:.4);
\draw (1,1.5)++(-20:.4) node[right] {1};
        \draw (1,1.5)--++(20:.4);
\draw (1,1.5)++(20:.4) node[right]
        {2};
\draw (1,1.5)--++(60:.4);
\draw (1,1.5)++(60:.4)
        node[right] {3};
      \end{tikzpicture}
    \\$(\gamma_{\mathrm{III}})$&$(\gamma_{\mathrm{IV}})$\\
    &\\
  \end{tabular}
  
  \caption{Strata that intersect $\overline{\Per_{2,5}}$}
  \label{fig:StrataThatIntersectC5}  
\end{table}

\begin{table}
    \centering
    \begin{tabular}{ccc}
    \begin{tikzpicture}
      \draw (-1.5,0) node {$\tau$};
\draw (-1.5,2) node {$\sigma$};
\draw
      (0,0)--(3,0);
\draw (0,0) node {$\bullet$};
\draw
      (0,0)--++(160:.4);
\draw (0,0)++(160:.4) node[left] {2};
\draw
      (0,0)--++(200:.4);
\draw (0,0)++(200:.4) node[left] {$*$};
\draw
      (1,0) node {$\bullet$};
\draw (1,0)--++(-90:.4);
\draw
      (1,0)++(-90:.4) node[below] {4};
\draw (2,0) node {$\bullet$};
      \draw (2,0)--++(-90:.4);
\draw (2,0)++(-90:.4) node[below] {3};
      \draw (3,0) node {$\bullet$};
\draw (3,0)--++(20:.4);
\draw
      (3,0)++(20:.4) node[right] {1};
\draw (3,0)--++(-20:.4);
\draw
      (3,0)++(-20:.4) node[right] {5};
\draw[->]
      (1.5,.75)--(1.5,0.25);
\draw
      (3,2.5)--(2,2.5)--(1,2.5)--(0,2)--(1,1.5)--(2,1.5)--(3,1.5);
      \draw (3,2.5) node {$\bullet$};
\draw (3,2.5)--++(20:.4);
\draw
      (3,2.5)++(20:.4) node[right] {5};
\draw (2,2.5) node
      {$\bullet$};
\draw (1,2.5) node {$\bullet$};
\draw (0,2) node
      {$\bullet$};
\draw[very thick] (0,2)--++(160:.4);
\draw
      (0,2)++(160:.4) node[left] {1};
\draw[very thick]
      (0,2)--++(200:.4);
\draw (0,2)++(200:.4) node[left] {$*$};
\draw
      (1,1.5) node {$\bullet$};
\draw (1,1.5)--++(-90:.4);
\draw
      (1,1.5)++(-90:.4) node[below] {3};
\draw (2,1.5) node
      {$\bullet$};
\draw (2,1.5)--++(-90:.4);
\draw (2,1.5)++(-90:.4)
      node[below] {2};
\draw (3,1.5) node {$\bullet$};
\draw
      (3,1.5)--++(-20:.4);
\draw (3,1.5)++(-20:.4) node[right] {4};
    \end{tikzpicture}&
    \begin{tikzpicture}
\draw
      (0,0)--(3,0);
\draw (0,0) node {$\bullet$};
\draw
      (0,0)--++(160:.4);
\draw (0,0)++(160:.4) node[left] {2};
\draw
      (0,0)--++(200:.4);
\draw (0,0)++(200:.4) node[left] {$*$};
\draw
      (1,0) node {$\bullet$};
\draw (1,0)--++(-90:.4);
\draw
      (1,0)++(-90:.4) node[below] {4};
\draw (2,0) node {$\bullet$};
      \draw (2,0)--++(-90:.4);
\draw (2,0)++(-90:.4) node[below] {1};
      \draw (3,0) node {$\bullet$};
\draw (3,0)--++(20:.4);
\draw
      (3,0)++(20:.4) node[right] {3};
\draw (3,0)--++(-20:.4);
\draw
      (3,0)++(-20:.4) node[right] {5};
\draw[->]
      (1.5,.75)--(1.5,0.25);
\draw
      (3,2.5)--(2,2.5)--(1,2.5)--(0,2)--(1,1.5)--(2,1.5)--(3,1.5);
      \draw (3,2.5) node {$\bullet$};
\draw (2,2.5) node
      {$\bullet$};
      \draw (2,2.5)--++(-90:.4);
\draw
      (2,2.5)++(-90:.4) node[below] {5};
\draw (1,2.5) node {$\bullet$};
\draw (1,2.5)--++(-90:.4);
\draw
      (1,2.5)++(-90:.4) node[below] {3};
\draw (0,2) node
      {$\bullet$};
\draw[very thick] (0,2)--++(160:.4);
\draw
      (0,2)++(160:.4) node[left] {1};
\draw[very thick]
      (0,2)--++(200:.4);
\draw (0,2)++(200:.4) node[left] {$*$};
\draw
      (1,1.5) node {$\bullet$};
\draw (2,1.5) node
      {$\bullet$};
\draw (3,1.5) node {$\bullet$};
\draw (3,1.5)--++(20:.4);
\draw (3,1.5)++(20:.4)
      node[right] {2};
\draw
      (3,1.5)--++(-20:.4);
\draw (3,1.5)++(-20:.4) node[right] {4};
    \end{tikzpicture}&
    \begin{tikzpicture}
\draw
      (0,0)--(3,0);
\draw (0,0) node {$\bullet$};
\draw
      (0,0)--++(160:.4);
\draw (0,0)++(160:.4) node[left] {2};
\draw
      (0,0)--++(200:.4);
\draw (0,0)++(200:.4) node[left] {4};
\draw
      (1,0) node {$\bullet$};
\draw (1,0)--++(-90:.4);
\draw
      (1,0)++(-90:.4) node[below] {$*$};
\draw (2,0) node {$\bullet$};
      \draw (2,0)--++(-90:.4);
\draw (2,0)++(-90:.4) node[below] {5};
      \draw (3,0) node {$\bullet$};
\draw (3,0)--++(20:.4);
\draw
      (3,0)++(20:.4) node[right] {1};
\draw (3,0)--++(-20:.4);
\draw
      (3,0)++(-20:.4) node[right] {3};
\draw[->]
      (1.5,.75)--(1.5,0.25);
\draw
      (3,2.5)--(2,2.5)--(1,2)--(0,2)--(1,2)--(2,1.5)--(3,1.5);
      \draw (3,2.5) node {$\bullet$};
\draw (3,2.5)--++(20:.4);
\draw
      (3,2.5)++(20:.4) node[right] {5};
\draw (2,2.5) node
      {$\bullet$};
\draw (0,2) node
      {$\bullet$};
\draw[very thick] (0,2)--++(160:.4);
\draw
      (0,2)++(160:.4) node[left] {1};
\draw
      (0,2)--++(200:.4);
\draw (0,2)++(200:.4) node[left] {3};
\draw[very thick] (0,2)--(1,2);
\draw
      (1,2) node {$\bullet$};
\draw[very thick] (1,2)--++(-90:.4);
\draw
      (1,2)++(-90:.4) node[below] {$*$};
\draw (2,1.5) node
      {$\bullet$};
\draw (2,1.5)--++(-90:.4);
\draw (2,1.5)++(-90:.4)
      node[below] {4};
\draw (3,1.5) node {$\bullet$};
\draw
      (3,1.5)--++(-20:.4);
\draw (3,1.5)++(-20:.4) node[right] {2};
    \end{tikzpicture}
    \\degen of $\gamma_4,\gamma_5$&degen of $\gamma_5$&degen of $\gamma_5$
    \\\hline
    \begin{tikzpicture}
          \draw (-1.5,0) node {$\tau$};
\draw (-1.5,2) node {$\sigma$};
\draw
      (0,0)--(2,0);
\draw (0,0) node {$\bullet$};
\draw
      (0,0)--++(160:.4);
\draw (0,0)++(160:.4) node[left] {4};
\draw
      (0,0)--++(200:.4);
\draw (0,0)++(200:.4) node[left] {$*$};
\draw
      (1,0) node {$\bullet$};
\draw (1,0)--++(-90:.4);
\draw
      (1,0)++(-90:.4) node[below] {2};
\draw (2,0) node {$\bullet$};
\draw (2,0)--++(-40:.4);
\draw (2,0)++(-40:.4) node[right] {5};
\draw (2,0)--++(0:.4);
\draw (2,0)++(0:.4) node[right] {3};
\draw (2,0)--++(40:.4);
\draw
      (2,0)++(40:.4) node[right] {1};
\draw[->]
      (1.5,.75)--(1.5,0.25);
\draw
      (2,2.5)--(1,2)--(0,2)--(1,2)--(2,1.5);
\draw (2,2.5) node
      {$\bullet$};
            \draw (2,2.5)--++(40:.4);
\draw
      (2,2.5)++(40:.4) node[right] {5};
\draw (1,2) node {$\bullet$};
\draw[very thick] (1,2)--++(-90:.4);
\draw (1,2)++(-90:.4) node[below] {1};
\draw[very thick] (0,2)--(1,2);
\draw (0,2) node
      {$\bullet$};
\draw (0,2)--++(160:.4);
\draw
      (0,2)++(160:.4) node[left] {3};
\draw[very thick]
      (0,2)--++(200:.4);
\draw (0,2)++(200:.4) node[left] {$*$};
\draw (2,1.5) node
      {$\bullet$};
\draw (2,1.5)--++(0:.4);
\draw (2,1.5)++(0:.4)
      node[right] {2};
\draw
      (2,1.5)--++(-40:.4);
\draw (2,1.5)++(-40:.4) node[right] {4};
    \end{tikzpicture}&
    \begin{tikzpicture}
\draw
      (0,0)--(3,0);
\draw (0,0) node {$\bullet$};
\draw
      (0,0)--++(160:.4);
\draw (0,0)++(160:.4) node[left] {2};
\draw
      (0,0)--++(200:.4);
\draw (0,0)++(200:.4) node[left] {$*$};
\draw
      (1,0) node {$\bullet$};
\draw (1,0)--++(-90:.4);
\draw
      (1,0)++(-90:.4) node[below] {5};
\draw (2,0) node {$\bullet$};
      \draw (2,0)--++(-90:.4);
\draw (2,0)++(-90:.4) node[below] {1};
      \draw (3,0) node {$\bullet$};
\draw (3,0)--++(20:.4);
\draw
      (3,0)++(20:.4) node[right] {3};
\draw (3,0)--++(-20:.4);
\draw
      (3,0)++(-20:.4) node[right] {4};
\draw[->]
      (1.5,.75)--(1.5,0.25);
\draw
      (3,2.5)--(2,2.5)--(1,2.5)--(0,2)--(1,1.5)--(2,1.5)--(3,1.5);
      \draw (3,2.5) node {$\bullet$};
      \draw (3,2.5)--++(20:.4);
\draw (3,2.5)++(20:.4)
      node[right] {2};
\draw (2,2.5) node
      {$\bullet$};
      \draw (2,2.5)--++(-90:.4);
\draw
      (2,2.5)++(-90:.4) node[below] {5};
\draw (1,2.5) node {$\bullet$};
\draw (0,2) node
      {$\bullet$};
\draw[very thick] (0,2)--++(160:.4);
\draw
      (0,2)++(160:.4) node[left] {1};
\draw[very thick]
      (0,2)--++(200:.4);
\draw (0,2)++(200:.4) node[left] {$*$};
\draw
      (1,1.5) node {$\bullet$};
\draw (1,1.5)--++(-90:.4);
\draw
      (1,1.5)++(-90:.4) node[below] {4};
\draw (2,1.5) node
      {$\bullet$};
\draw (3,1.5) node {$\bullet$};
\draw
      (3,1.5)--++(-20:.4);
\draw (3,1.5)++(-20:.4) node[right] {3};
    \end{tikzpicture}&
    \begin{tikzpicture}
\draw
      (0,0)--(3,0);
\draw (0,0) node {$\bullet$};
\draw
      (0,0)--++(160:.4);
\draw (0,0)++(160:.4) node[left] {2};
\draw
      (0,0)--++(200:.4);
\draw (0,0)++(200:.4) node[left] {$*$};
\draw
      (1,0) node {$\bullet$};
\draw (1,0)--++(-90:.4);
\draw
      (1,0)++(-90:.4) node[below] {5};
\draw (2,0) node {$\bullet$};
      \draw (2,0)--++(-90:.4);
\draw (2,0)++(-90:.4) node[below] {4};
      \draw (3,0) node {$\bullet$};
\draw (3,0)--++(20:.4);
\draw
      (3,0)++(20:.4) node[right] {1};
\draw (3,0)--++(-20:.4);
\draw
      (3,0)++(-20:.4) node[right] {3};
\draw[->]
      (1.5,.75)--(1.5,0.25);
\draw
      (3,2.5)--(2,2.5)--(1,2.5)--(0,2)--(1,1.5)--(2,1.5)--(3,1.5);
      \draw (3,2.5) node {$\bullet$};
\draw (3,2.5)--++(20:.4);
\draw
      (3,2.5)++(20:.4) node[right] {5};
      \draw (3,2.5)--++(-20:.4);
\draw
      (3,2.5)++(-20:.4) node[right] {2};
\draw (2,2.5) node
      {$\bullet$};
\draw (1,2.5) node {$\bullet$};
\draw (1,2.5)--++(-90:.4);
\draw
      (1,2.5)++(-90:.4) node[right] {4};
\draw (0,2) node
      {$\bullet$};
\draw[very thick] (0,2)--++(160:.4);
\draw
      (0,2)++(160:.4) node[left] {1};
\draw[very thick]
      (0,2)--++(200:.4);
\draw (0,2)++(200:.4) node[left] {$*$};
\draw
      (1,1.5) node {$\bullet$};
\draw (2,1.5) node
      {$\bullet$};
\draw (2,1.5)--++(-90:.4);
\draw (2,1.5)++(-90:.4)
      node[below] {3};
\draw (3,1.5) node {$\bullet$};
    \end{tikzpicture}
    \\degen of $\gamma_5,\gamma_{\mathrm{III}}$&degen of $\gamma_6$&degen of $\gamma_6$
    \\\hline
    \begin{tikzpicture}
          \draw (-1.5,0) node {$\tau$};
\draw (-1.5,2) node {$\sigma$};
\draw
      (0,0)--(3,0);
\draw (0,0) node {$\bullet$};
\draw
      (0,0)--++(160:.4);
\draw (0,0)++(160:.4) node[left] {2};
\draw
      (0,0)--++(200:.4);
\draw (0,0)++(200:.4) node[left] {5};
\draw
      (1,0) node {$\bullet$};
\draw (1,0)--++(-90:.4);
\draw
      (1,0)++(-90:.4) node[below] {$*$};
\draw (2,0) node {$\bullet$};
      \draw (2,0)--++(-90:.4);
\draw (2,0)++(-90:.4) node[below] {3};
      \draw (3,0) node {$\bullet$};
\draw (3,0)--++(20:.4);
\draw
      (3,0)++(20:.4) node[right] {1};
\draw (3,0)--++(-20:.4);
\draw
      (3,0)++(-20:.4) node[right] {4};
\draw[->]
      (1.5,.75)--(1.5,0.25);
\draw
      (3,2.5)--(2,2.5)--(1,2)--(2,1.5)--(3,1.5);
      \draw[very thick] (0,2)--(1,2);
      \draw (3,2.5) node {$\bullet$};
\draw (3,2.5)--++(20:.4);
\draw
      (3,2.5)++(20:.4) node[right] {5};
\draw (2,2.5) node
      {$\bullet$};
      \draw (2,2.5)--++(-90:.4);
\draw (2,2.5)++(-90:.4)
      node[below] {2};
\draw (0,2) node
      {$\bullet$};
\draw[very thick] (0,2)--++(160:.4);
\draw
      (0,2)++(160:.4) node[left] {1};
\draw
      (0,2)--++(200:.4);
\draw (0,2)++(200:.4) node[left] {4};
\draw
      (1,2) node {$\bullet$};
\draw[very thick] (1,2)--++(-90:.4);
\draw
      (1,2)++(-90:.4) node[below] {$*$};
\draw (2,1.5) node
      {$\bullet$};
\draw (3,1.5) node {$\bullet$};
\draw
      (3,1.5)--++(-20:.4);
\draw (3,1.5)++(-20:.4) node[right] {3};
    \end{tikzpicture}&
    \begin{tikzpicture}
\draw
      (0,0)--(2,0);
\draw (0,0) node {$\bullet$};
\draw
      (0,0)--++(160:.4);
\draw (0,0)++(160:.4) node[left] {5};
\draw
      (0,0)--++(200:.4);
\draw (0,0)++(200:.4) node[left] {$*$};
\draw
      (1,0) node {$\bullet$};
\draw (1,0)--++(-90:.4);
\draw
      (1,0)++(-90:.4) node[below] {2};
\draw (2,0) node {$\bullet$};
\draw (2,0)--++(-40:.4);
\draw (2,0)++(-40:.4) node[right] {4};
\draw (2,0)--++(0:.4);
\draw (2,0)++(0:.4) node[right] {3};
\draw (2,0)--++(40:.4);
\draw
      (2,0)++(40:.4) node[right] {1};
\draw[->]
      (1.5,.75)--(1.5,0.25);
\draw
      (2,2.5)--(1,2)--(2,1.5);
      \draw[very thick] (0,2)--(1,2);
\draw (2,2.5) node
      {$\bullet$};
            \draw (2,2.5)--++(40:.4);
\draw
      (2,2.5)++(40:.4) node[right] {5};
      \draw (2,2.5)--++(0:.4);
\draw (2,2.5)++(0:.4)
      node[right] {2};
\draw (1,2) node {$\bullet$};
\draw[very thick] (1,2)--++(-90:.4);
\draw
      (1,2)++(-90:.4) node[below] {1};
\draw (0,2) node
      {$\bullet$};
\draw (0,2)--++(160:.4);
\draw
      (0,2)++(160:.4) node[left] {4};
\draw[very thick]
      (0,2)--++(200:.4);
\draw (0,2)++(200:.4) node[left] {$*$};
\draw (2,1.5) node
      {$\bullet$};

\draw
      (2,1.5)--++(-40:.4);
\draw (2,1.5)++(-40:.4) node[right] {3};
    \end{tikzpicture}&
    \fbox{
    \begin{tikzpicture}
        \draw (0,0)--(1,0);
\draw (0,0) node
        {$\bullet$};
\draw (0,0)--++(140:.4);
\draw (0,0)++(140:.4)
        node[left] {$*$};
\draw (0,0)--++(180:.4);
\draw (0,0)++(180:.4)
        node[left] {1};
\draw (0,0)--++(220:.4);
\draw (0,0)++(220:.4)
        node[left] {3};
\draw (1,0) node
        {$\bullet$};
\draw (1,0)--++(-40:.4);
\draw (1,0)++(-40:.4)
        node[right] {2};
\draw (1,0)--++(0:.4);
\draw
        (1,0)++(0:.4) node[right] {4};
\draw (1,0)--++(40:.4);
\draw
        (1,0)++(40:.4) node[right] {5};

\draw[->] (.5,1)--(.5,0.5);
        \draw[very thick] (1,1.5)--(0,1.5);
\draw (0,1.5) node
        {$\bullet$};
\draw (0,1.5)--++(140:.4);
\draw
        (0,1.5)++(140:.4) node[left] {$*$};
        \draw (0,1.5)--++(180:.4);
\draw
        (0,1.5)++(180:.4) node[left] {5};
\draw[very thick]
        (0,1.5)--++(220:.4);
\draw (0,1.5)++(220:.4) node[left]
        {2};
\draw (1,1.5) node
        {$\bullet$};
\draw (1,1.5)--++(-40:.4);
\draw
        (1,1.5)++(-40:.4) node[right] {1};
\draw[very thick]
        (1,1.5)--++(0:.4);
\draw (1,1.5)++(0:.4) node[right] {3};
        \draw (1,1.5)--++(40:.4);
\draw (1,1.5)++(40:.4) node[right]
        {4};
      \end{tikzpicture}
      }
    \\degen of $\gamma_6$&degen of $\gamma_6,\gamma_{\mathrm{IV}}$&
    \\
    &\\
    \end{tabular}
    \caption{Strata that intersect $\tilde{\Delta}_{2,5}$, but do not intersect $\overline{\Per_{2,5}}$ (except the last, see Section \ref{sec:R6})}
    \label{fig:DoNotIntersectC5}
\end{table}

\subsection{Local analysis of \texorpdfstring{$\tilde{\Delta}_{2,5}$}{Delta\_\{2,5\}-tilde} at boundary points}\label{sec:LocalAnalysis} For convenience, we write $R_1,\ldots,R_{\mathrm{IV}}$ and $U_1,\ldots,U_{\mathrm{IV}}$ for $R_{\gamma_1},\ldots,R_{\gamma_{\mathrm{IV}}}$ and $U_{\gamma_1},\ldots,U_{\gamma_{\mathrm{IV}}}$, respectively, where $U_\gamma$ is the infinitesimal neighborhood of $R_\gamma\subseteq\Hbar_{2,5}.$

\subsubsection{Local analysis near $R_1$} 
Let $s_1,$ $s_2,$ and $s_3$ be node smoothing parameters for edges of $\sigma$ as indicated in Table \ref{fig:StrataThatIntersectC5}. By Section \ref{sec:HnBarCoords}, these are local coordinates on the infinitesimal neighborhood $U_1$ of $R_1\subseteq\Hbar_{2,5}.$ By Section \ref{sec:pisanddeltacoords},  $\tilde{\Delta}_{2,5}\subseteq\Hbar_{2,5}$ is defined (on $U_1$) by $\CR(C,a_1,a_{3,0},a_{5,0},a_{4,0})-\CR(D,b_1,b_3,b_5,b_4)$ and $\CR(C,a_1,a_{3,0},a_{2,0},a_{5,0})-\CR(D,b_1,b_3,b_2,b_5).$ By Observations \ref{obs:DifferentChoicesOfSmoothingParameter} and \ref{obs:SmoothingParametersOnSourceAndTarget}, we have
\begin{align*}
    \CR(C,a_1,a_{3,0},a_{5,0},a_{4,0})&\sim s_2,&\CR(D,b_1,b_3,b_5,b_4)&\sim s_3,\\\CR(C,a_1,a_{3,0},a_{2,0},a_{5,0})&\sim s_1,&\CR(D,b_1,b_3,b_2,b_5)&\sim s_2.
\end{align*}
Thus $\tilde{\Delta}_{2,5}$ is locally defined by $s_2-\alpha s_3$ and $s_1-\beta s_2$ for nonvanishing functions $\alpha,\beta$ on $R_1.$ It follows that $p_1:=R_1$ is a smooth point of $\overline{\Per_{2,5}},$ and that at $p_1,$ $\overline{\Per_{2,5}}$ is transverse to the three boundary hypersurfaces defined by $s_1,$ $s_2,$ and $s_3$, respectively.

\subsubsection{Local analysis near $R_2$} Let $s_1,$ $s_2,$ and $s_3$ be node smoothing parameters for edges of $\sigma$ as indicated in Table \ref{fig:StrataThatIntersectC5}; these are coordinates on $U_2$. Locally, $\tilde{\Delta}_{2,5}$ is defined by $\CR(C,a_{2,0},a_1,a_{3,0},a_{5,0})-\CR(D,b_2,b_1,b_3,b_5)$ and $\CR(C,a_1,a_{5,0},a_{3,0},a_{4,0})-\CR(D,b_1,b_5,b_3,b_4).$ By Observations \ref{obs:DifferentChoicesOfSmoothingParameter}, \ref{obs:ProductOfSmoothingParameters}, and \ref{obs:SmoothingParametersOnSourceAndTarget}, we have
\begin{align*}
    \CR(C,a_{2,0},a_1,a_{3,0},a_{5,0})&\sim s_1,&\CR(D,b_2,b_1,b_3,b_5)&\sim s_2,\\\CR(C,a_1,a_{5,0},a_{3,0},a_{4,0})&\sim s_1s_2,&\CR(D,b_1,b_5,b_3,b_4)&\sim s_3.
\end{align*}
Thus $\tilde{\Delta}_{2,5}$ is locally defined by $s_1-\alpha s_2$ and $s_1s_2-\beta s_3$ for nonvanishing functions $\alpha,\beta$ on $R_2.$ It follows that $p_2:=R_2$ is a smooth point of $\overline{\Per_{2,5}},$ and that at $p_2$, $\overline{\Per_{2,5}}$ is transverse to the hypersurfaces defined by $s_1$ and $s_2,$ and is (simply) tangent to the boundary hypersurface  defined by $s_3$.

\subsubsection{Local analysis near $R_3$} The analysis is identical to that of $R_1$; $R_3$ is a smooth point of $\overline{\Per_{2,5}}$ at which $\overline{\Per_{2,5}}$ is transverse to the three boundary hypersurfaces of $\Hbar_{2,5}$ that intersect at $R_1$.

\subsubsection{Local analysis near $R_4$} Let $s_1$ and $s_2$ be node smoothing parameters for edges of $\sigma$ as indicated in Table \ref{fig:StrataThatIntersectC5}. Let $s_3=\CR(C,a_1,a_{2,0},a_{3,0},a_{4,0})$. (This is a cross-ratio parameter for the lower middle vertex of $\sigma$.) Then $s_1,s_2,s_3$ are coordinates for $U_4$. Note that on $U_4,$ we have $R_4=Z(s_1,s_2)$. Locally, $\tilde{\Delta}_{2,5}$ is defined by $\CR(C,a_{2,0},a_1,a_{3,0},a_{5,0})-\CR(D,b_2,b_1,b_3,b_5)$ and $\CR(C,a_1,a_{2,0},a_{3,0},a_{4,0})-\CR(D,b_1,b_2,b_3,b_4).$ By Observations \ref{obs:DifferentChoicesOfSmoothingParameter} and \ref{obs:SmoothingParametersOnSourceAndTarget}, we have $\CR(C,a_{2,0},a_1,a_{3,0},a_{5,0})\sim s_1$ and $\CR(D,b_2,b_1,b_3,b_5)\sim s_2.$
By Observation \ref{obs:CrossRatiosParametersOnSourceAndTarget}, the function $g:=\CR(D,b_1,b_2,b_3,b_4)-\CR(C,a_{4,0},a_1,a_{2,0},a_{3,0})$ on $U_4$ is in the ideal $(s_1,s_2)$. On the other hand, the cross-ratio functional equations give $$\CR(C,a_4,a_1,a_2,a_3)=\frac{\CR(C,a_1,a_2,a_3,a_4)}{\CR(C,a_1,a_2,a_3,a_4)-1}=\frac{s_3}{s_3-1}.$$ Thus
$$\CR(C,a_{4,0},a_1,a_{2,0},a_{3,0})-\CR(D,b_1,b_2,b_3,b_4)=s_3-\frac{s_3}{s_3-1}-g=\frac{s_3(s_3-2)}{s_3-1}-g$$
Thus $\tilde{\Delta}_{2,5}$ is locally defined by $s_1-\alpha s_2$ and $\frac{s_3(s_3-2)}{s_3-1}-g$, where $\alpha$ is nonvanishing on $R_4$ and $g\in(s_1,s_2).$ Since $s_3\ne0$ on $R_4,$ $\overline{\Per_{2,5}}\cap R_4$ is a single reduced point $p_4$ at which $s_3=\CR(C,a_1,a_{2,0},a_{3,0},a_{4,0})=2$. 

In fact, reducedness in this case implies that $p_4$ is a smooth point of $\overline{\Per_{2,5}}$ at which $\overline{\Per_{2,5}}$ is transverse to the boundary hypersurfaces defined by $s_1$ and $s_2,$ respectively. This can be seen either by a standard Jacobian computation, or by substituting $s_1\mapsto\alpha s_2$ so that $\tilde{\Delta}_{2,5}$ and $R_4$ are curves inside a smooth surface.

\subsubsection{Local analysis near $R_5$} 
We have coordinates $s_1=\CR(C,a_1,a_{4,0},a_{3,0},a_{2,0}),$ $s_2=\CR(C,a_1,a_*,a_{3,0},a_{4,0})$, and $s_3=\CR(C,a_*,a_{4,0},a_{5,1},a_,{2,0})$ on $U_5$. ) Here $s_1$ is a node smoothing parameter, and $s_2$ and $s_3$ are cross-ratio parameters on the left and lower-right vertices of $\sigma$, respectively; thus $R_5$ is defined (on $U_5$) by the vanishing of $s_1$. (Note: This is the first time that we are forced to use the enlarged marking set of Remark \ref{rem:BaseChange}; this is due to the fact the no vertex on the right side of $\sigma$ has valence 4.) One may check that the map $(s_1,s_2,s_3):\H_{2,5}\to\A^3$ is generically injective, hence Caution \ref{cau:Roots} does not apply.

Locally, $\tilde{\Delta}_{2,5}$ is defined by
\begin{align*}
    h_1:&=\CR(C,a_1,a_{4,0},a_{3,0},a_{2,0})-\CR(D,b_1,b_4,b_3,b_2)\quad\quad\text{and}\\
    h_2:&=\CR(C,a_1,a_{3,0},a_{5,0},a_{4,0})-\CR(D,b_1,b_3,b_5,b_4).
\end{align*}
The cross-ratio functional equations and the relations \eqref{eq:RelsBtwnCR} imply the following coordinate expressions\footnote{Observation \ref{obs:SmoothingParametersOnSourceAndTarget} implies that $\CR(D,b_1,b_4,b_3,b_2)\sim s_1$, but this is not helpful here! This is because analyzing the equation $h_1=s_1-\alpha s_1$ requires specific information about the nonvanishing function $\alpha$. This occurs whenever vertices $v\in\sigma$ $\phi(v)\in\tau$ have ``the same'' node smoothing parameter --- it is also why we chose a \textit{specific} node smoothing parameter.} on $U_5$:
\begin{align*}
    h_1&=s_1\left(\frac{-2s_2+s_3+3s_2s_3}{s_3(1+s_2)}\right)+s_1g_1\\
    h_2&=\frac{-1-s_2+s_3-s_2s_3}{(1+s_2)s_3}+g_2,
\end{align*}
where $g_1,g_2\in(s_1)$.
Note that $Z(h_1)\supseteq R_5$, and $Z(h_2)\cap R_5=Z(h_2,s_1)$ is a curve; in particular, $\tilde{\Delta}_{2,5}\cap R_5$ is the curve $Z(h_2)\cap R_5$. On the other hand, since $h_1$ vanishes to order 1 along $R_5$, $h_1$ defines the hypersurface $R_5\cup Z(h_1/s_1)$ in $U_5$. 
By Thurston/Epstein transversality, $Z(h_2,h_1/s_1)$ is a curve containing $\overline{\Per_{2,5}}$ (and possibly components contained in $R_5$). 
We compute that $R_4\cap\overline{\Per_{2,5}}=R_5\cap Z(h_2)\cap Z(h_1/s_1)=Z(s_1,h_2,h_1/s_1)$ consists of the two reduced points $p_5,p_5'$ given by $(s_1,s_2,s_3)=(0,\frac{-1\pm2i}{5},\frac{1\pm i}{2})$. Since $R_5$ is a hypersurface, reducedness implies that these are smooth points of $\overline{\Per_{2,5}}$ at which $\overline{\Per_{2,5}}$ and $R_5$ intersect transversely. For convenience later we record the cross-ratios $\CR(C,a_1,a_{2,0},a_{3,0},a_{5,0})=\frac{1\pm i}{2}$ and $\CR(C,a_1,a_{3,0},a_{4,0},a_{5,0})=\frac{1\mp i}{2}$.

\subsubsection{Local analysis near $R_6$}\label{sec:R6} We have coordinates $s_1=\CR(C,a_1,a_{5,0},a_{4,0},a_{2,0})$, $s_2=\CR(C,a_1,a_*,a_{4,0},a_{3,0}),$ and $s_3=\CR(C,a_*,a_{2,0},a_{3,1},a_{5,0})$ on $U_6$. 
Note that $s_1$ is a node smoothing parameter, and $s_2$ and $s_3$ are cross-ratio parameters on the left and upper-right vertices of $\sigma$, respectively; thus $R_6$ is defined (on $U_6$) by the vanishing of $s_1$. One may check that the map $(s_1,s_2,s_3):\Hbar_{2,5}\to\A^3$ is generically 2-to-1, hence Caution \ref{cau:Roots} \emph{does} apply\footnote{It happens that there do exist Zariski coordinates around $R_6$, making the analysis similar to that of $R_5$. We chose these coordinates instead to demonstrate the most general form of the process, and why it is necessary to work with power series.}. While this map is guaranteed to restrict to an isomorphism on $U_6$, it is indeed ramified over the locus $Z(q)\subseteq\A^3$, where $$q=4 s_1 s_2 (-1 + s_3) s_3 + (s_2 s_3 - s_1 (s_2 + s_3))^2.$$ Correspondingly, there is a second boundary stratum $R_{\gamma}$ such that $s_1,s_2,s_3$ are also coordinates on $U_{\gamma}$ --- see the boxed stratum in Table \ref{fig:DoNotIntersectC5}. 

Locally, $\tilde{\Delta}_{2,5}$ is defined by 
\begin{align*}
    h_1:&=\CR(C,a_1,a_{5,0},a_{4,0},a_{2,0})-\CR(D,b_1,b_5,b_4,b_2)\quad\quad\text{and}\\
    h_2:&=\CR(C,a_1,a_{3,0},a_{4,0},a_{5,0})-\CR(D,b_1,b_3,b_4,b_5).
\end{align*} Naively applying the cross-ratio functional equations to write $h_1$ and $h_2$ in terms of $s_1,s_2,s_3$ yields expressions involving $\sqrt{q}$. Writing $h_1$ and $h_2$ as elements of $k[s_2,s_3][[s_1]]$ in terms of the two choices of $\sqrt{q}$ give their local expansions along $R_6$ and $R_\gamma$, respectively. For example, the choices of $\sqrt{q}$ give rise to the two expressions 
\begin{align*}
    h_1&=\frac{3s_2s_3-2 s_2-s_3}{(s_2-1) s_3}s_1+s_1g_1&&\text{or}&
    h_1'&=-\frac{s_2^2}{s_2^2-1}+s_1+s_1g_2,
\end{align*}
for $g_1,g_2\in(s_1)$. Note that the first expression is the correct one (hence the naming choice), since we know that $h_1$ vanishes along $R_6$. (One could also have distinguished between the two square roots in other ways, e.g. by noting that on $R_6$ we must have $\CR(C,a_*,a_1,a_{2,0},a_{3,0})=-1.$) We apply similar reasoning to $h_2$ to find the correct expansion $$h_2=\frac{1+s_2-s_3+s_2s_3}{-1+s_2}+s_1g_3,$$ where $g_3\in(s_1).$
By the same reasoning as in the analysis of $R_5,$ $\overline{\Per_{2,5}}\cap R_6=Z(s_1,h_1/s_1,h_2)$; this consists of the two reduced points $p_6,p_6'$ given by $(s_1,s_2,s_3)=(0,\frac{\pm1}{\sqrt{5}},\frac{3\pm\sqrt{5}}{2})$. These are smooth points of $\overline{\Per_{2,5}}$ at which $\overline{\Per_{2,5}}$ meets $R_6$ transversely. For convenience we record the cross-ratios $\CR(C,a_1,a_{3,0},a_{4,0},a_{5,0})=\frac{-1\pm\sqrt{5}}{2}$.

\subsubsection{Local analysis near $R_7$} Note that $\pi_1$ and $\pi_2$ map $R_7$ to $M_{0,5}.$ Let $s_2=\CR(C,a_1,a_{2,0},a_{3,0},a_{4,0})$ and $s_3=\CR(C,a_1,a_{2,0},a_{3,0},a_{5,0})$. Then $s_2$ and $s_3$ are global coordinates on $R_7$. Locally, $\tilde{\Delta}_{2,5}$ is defined by $\CR(C,a_1,a_{2,0},a_{3,0},a_{4,0})-\CR(D,b_1,b_2,b_3,b_4)$ and $\CR(C,a_1,a_{2,0},a_{3,0},a_{5,0})-\CR(D,b_1,b_2,b_3,b_5).$ By Observation \ref{obs:CrossRatiosParametersOnSourceAndTarget} and the cross-ratio functional equations, on $R_7$ we have
\begin{align*}
    \CR(D,b_1,b_2,b_3,b_4)&=\CR(C,a_{5,0},a_1,a_{2,0},a_{3,0})=\frac{s_3}{s_3-1}\\
    \CR(D,b_1,b_2,b_3,b_5)&=\CR(C,a_{5,0},a_1,a_{2,0},a_{4,0})=\frac{s_3}{s_3-s_2}.
\end{align*}
Solving $h_1=h_2=0$ on $R_7$ gives the two reduced points $p_7,p_7'$ with $(s_2,s_3)=(\frac{1\pm\sqrt{5}}{2},\frac{3\pm\sqrt{5}}{2})$. As above, $p_7$ and $p_7'$ are smooth points of $\overline{\Per_{2,5}}$ at which $\overline{\Per_{2,5}}$ is transverse to $R_7$.

\subsubsection{Local analysis near $R_{\mathrm{I}},R_{\mathrm{II}},R_{\mathrm{III}},R_{\mathrm{IV}}$} These analyses are similar and the answers are identical. We work out $R_{\mathrm{I}}$ fully.

As with $R_7$, $\pi_1$ and $\pi_2$ map $R_{\mathrm{I}}$ to $M_{0,5}.$ Let $s_2=\CR(C,a_{5,0},a_1,a_{2,0},a_{3,0})$ and $s_3=\CR(a_{5,0},a_1,a_{2,0},a_{4,0})$, and let $s_1$ be any node smoothing parameter for the unique edge of $\sigma$. Then $s_1,s_2, s_3$ are coordinates on $U_{\mathrm{I}}$, with $s_2$ and $s_3$ restricting to global coordinates on $R_{\mathrm{I}}$. By Observation \ref{obs:CrossRatiosParametersOnSourceAndTarget} and the cross-ratio functional equations, restricted to $R_{\mathrm{I}}$ we have 
\begin{align*}
    \CR(D,b_5,b_1,b_2,b_3)=\CR(C,a_1,a_{2,0},a_{3,0},a_{4,0})^2&=\left(\frac{s_2(s_3-1)}{s_3(s_2-1)}\right)^2,
    \\
    \CR(D,b_5,b_1,b_2,b_4)=\CR(C,a_1,a_{2,0},a_{3,0},a_{5,0})^2&=\left(\frac{s_2}{s_2-1}\right)^2
    \\
\end{align*}
This means that, on $U_{\mathrm{I}}$, we have
\begin{align*}
    h_1:=\CR(C,a_{5,0},a_1,a_{2,0},a_{3,0})-  \CR(D,b_5,b_1,b_2,b_3)&=s_2-\left(\frac{s_2(s_3-1)}{s_3(s_2-1)}\right)^2+s_1g_1\\
h_2:=\CR(a_{5,0},a_1,a_{2,0},a_{4,0})-\CR(D,b_5,b_1,b_2,b_4)&=s_3-\left(\frac{s_2}{s_2-1}\right)^2+s_1g_2,
\end{align*}
where $g_1,g_2\in(s_1).$ We compute that $Z(h_1,h_2,s_1)=\tilde{\Delta}_{2,5}\cap R_{\mathrm{I}}$ consists of 5 reduced points. Again, since $R_{\mathrm{I}}$ is a hypersurface, these are smooth points of $\overline{\Per_{2,5}}$ at which $\overline{\Per_{2,5}}$ is transverse to $R_{\mathrm{I}}$.


\subsubsection{Local analysis near strata in Table \ref{fig:DoNotIntersectC5} (other than the boxed stratum).}
The strata in Table \ref{fig:DoNotIntersectC5} do not intersect $\overline{\Per_{2,5}}$; instead, they intersect the extra components of $\tilde{\Delta}_{2,5}$ found in the analyses of $R_5$ and $R_6$ above. We could show this directly for each stratum. Instead we consider the composition of $\pi_1$ with the forgetful map $\Mbar_{0,5}\to\Mbar_{0,4}$ that forgets the 5th marked point. This gives a map $\rho:\overline{\Per_{2,5}}\to\Mbar_{0,4}\cong\P^1$ of smooth curves. 

Using our complete enumeration of boundary strata in $\Hbar_{2,5},$ we easily check that $\rho^{-1}(12|34)$ is the single point $p_1.$ The node-smoothing parameter $\CR(1,2,3,4)$ at $12|34$ pulls back to $s_1s_2\sim s_1^2$ in the local coordinates given in the analysis of $R_1$. As $s_1$ is a local coordinate on $\overline{\Per_{2,5}}$, we conclude that $\rho$ is simply ramified at $p_1.$ We conclude that $\rho$ is a degree-2 map. We already know $\{p_5,p_5'\}\subseteq\rho^{-1}(13|24),$ so in fact $\{p_5,p_5'\}=\rho^{-1}(13|24).$ This implies that the degenerations of $\gamma_5$ in Table \ref{fig:DoNotIntersectC5} do not intersect $\overline{\Per_{2,5}}.$

An identical argument involving the forgetful map that forgets the 4th marked point shows that the degenerations of $\gamma_6$ in Table \ref{fig:DoNotIntersectC5} do not intersect $\overline{\Per_{2,5}}.$

\medskip

Summarizing this section, we have:
\begin{thm}
$\overline{\Per_{2,5}}$ is a smooth curve in $\Hbar_{2,5}$. The intersection of $\overline{\Per_{2,5}}$ with the boundary $\Hbar_{2,5}\setminus\H_{2,5}$ consists of the 10 points $p_1,p_2,p_3,p_4,p_5,p_5',p_6,p_6',p_7,p_7'$ --- the punctures of $\Per_{2,5}$ --- as well as the 20 points where $\overline{\Per_{2,5}}$ intersects $R_{\mathrm{I}},\ldots,R_{\mathrm{IV}}$ --- the points of $\Per_{2,5}\setminus\Per_{2,5}^\circ$, corresponding to maps $f:\P^1\to\P^1$ such that the unmarked critical point is in the orbit of the marked 5-periodic critical point.
\end{thm}

\subsection{\texorpdfstring{$\overline{\Per_{2,5}}$}{Per\_\{2,5\}-bar} as a cubic curve in \texorpdfstring{$\P^2$}{P2}}\label{sec:CubicCurve}
\begin{figure}
    \centering
    \begin{tikzpicture}[scale=0.9]
      \draw[thick,blue] (-30:.5) to[out=-90,in=45] (-90:3) to[out=-135,in=0, edge node={node[above] {$125|34$}}] (-135:4);
      \draw (90:.5) to[out=30,in=165] (30:3) to[out=-15,in=120, edge node={node[left] {$135|24$}}] (-15:4);
      \draw[thick,blue] (210:.5) to[out=150,in=285] (150:3) to[out=105,in=240, edge node={node[below right] {$145|23$}}] (105:4);
      \draw (-125:4) to[out=90,in=0] (-160:4) to[out=180,in=-90, edge node={node[above left] {$12|345$}}] (110:6);
      \draw[thick,blue] (-5:4) to[out=210,in=120] (-40:4) to[out=300,in=30, edge node={node[below] {$13|245$}}] (230:6);
      \draw (115:4) to[out=330,in=240] (80:4) to[out=420,in=150, edge node={node[right] {$14|235$}}] (350:6);
      \draw (-210:2) to[out=-180,in=90] (-160:4) to[out=-90,in=90, edge node={node[left] {$123|45$}}] (-120:6);
      \draw (-90:2) to[out=-60,in=210] (-40:4) to[out=30,in=210, edge node={node[below right] {$134|25$}}] (0:6);
      \draw[thick,blue] (30:2) to[out=60,in=330] (80:4) to[out=150,in=330, edge node={node[above right] {$124|35$}}] (120:6);
      \draw (0:1) to[out=90,in=0] (90:1) node[above=4pt,left=-3pt] {$15|234$} to[out=180,in=90] (180:1) to[out=270,in=180] (270:1) to[out=0,in=270] (0:1);
      \draw (-129:3.72) node {$\circ$} node[above=7pt,left=-2pt] {$p_1$};
      \draw[red,very thick] (-129:4.32)--(-129:3.12);
      \draw (159:2.48) node {$\circ$} node[below] {$p_2$};
      \draw[red,very thick] (159:3.08) to[out=-30,in=195] (159:2.48) to[out=15,in=260] (145:2.38);
      \draw (80:4) node {$\circ$} node[below=2pt,left=0pt] {$p_3$};
      \draw[red,very thick] (80:4)--++(-80:.6);
      \draw[red,very thick] (80:4)--++(100:.6);
      \draw (-100:1) node {$\circ$} node[below=2pt,left=0pt] {$p_4$};
      \draw[red,very thick] (-100:.5)--(-100:1.5);
      \draw (20:3.33) node {$\circ$} node[above=2pt,left=0pt] {$p_5$};
      \draw[red,very thick] (20:3.33)--++(50:.5);
      \draw[red,very thick] (20:3.33)--++(230:.5);
      \draw (30:3) node {$\circ$} node[above=8pt,left=0pt] {$p_5'$};
      \draw[red,very thick] (30:3)--++(70:.5);
      \draw[red,very thick] (30:3)--++(250:.5);
      \draw (-32:4.21) node {$\circ$} node[below=2pt] {$p_6$};
      \draw[red,very thick] (-32:4.21)--++(130:.5);
      \draw[red,very thick] (-32:4.21)--++(-50:.5);
      \draw (-25:4.41) node {$\circ$} node[above=2pt] {$p_6'$};
      \draw[red,very thick] (-25:4.41)--++(150:.5);
      \draw[red,very thick] (-25:4.41)--++(-30:.5);
    \end{tikzpicture}
    
    \caption{Diagram of how $\pi_1(\overline{\Per_{2,5}})\subseteq\Mbar_{0,5}$ intersects boundary strata}
    \label{fig:PictureOfCurve}
  \end{figure}
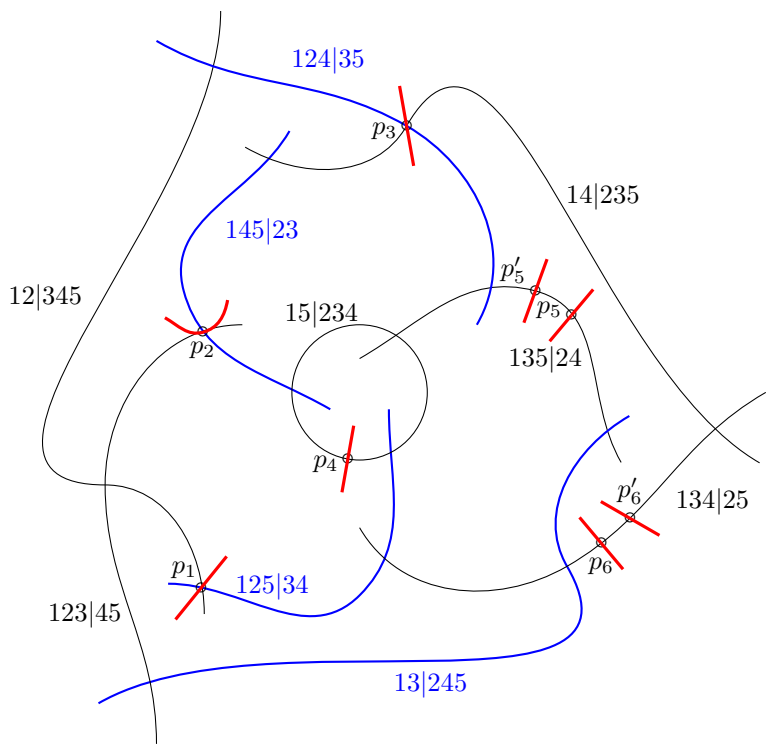

We now consider the image $\pi_1(\overline{\Per_{2,5}})=\pi_2(\overline{\Per_{2,5}})\subseteq\Mbar_{0,5}.$ Altogether, the analyses from the last section give the intersection points of $\overline{\Per_{2,5}}$ with the boundary of $\Mbar_{0,5}$ shown (with tangency drawn) in Figure \ref{fig:PictureOfCurve}. Here we use the abbreviation e.g. $12|345$ for the one-edge 5-marked tree (or corresponding closed boundary stratum) whose vertices are marked with $\{1,2\}$ and $\{3,4,5\},$ respectively. It is well-known \cite{Kapranov1993} that $\Mbar_{0,5}$ admits blow-down maps to $\P^2$; one such map, denoted $\psi_3$, contracts the four non-intersecting closed boundary strata $124|35,$ $125|34$, $145|23,$ and $13|245$ to points. (These curves are shown in blue in Figure \ref{fig:PictureOfCurve}.)

  Note that each point of $\overline{\Per_{2,5}}$ on a blue curve $\bar{Q_{\sigma}}$ is at an intersection point of $\bar{Q_{\sigma}}$ with another closed boundary curve $\bar{Q_{\sigma'}}$. Blowing down $\bar{Q_{\sigma}}$ has the effect of increasing the order of contact between $\overline{\Per_{2,5}}$ and $\bar{Q_{\sigma'}}$ by 1.
  
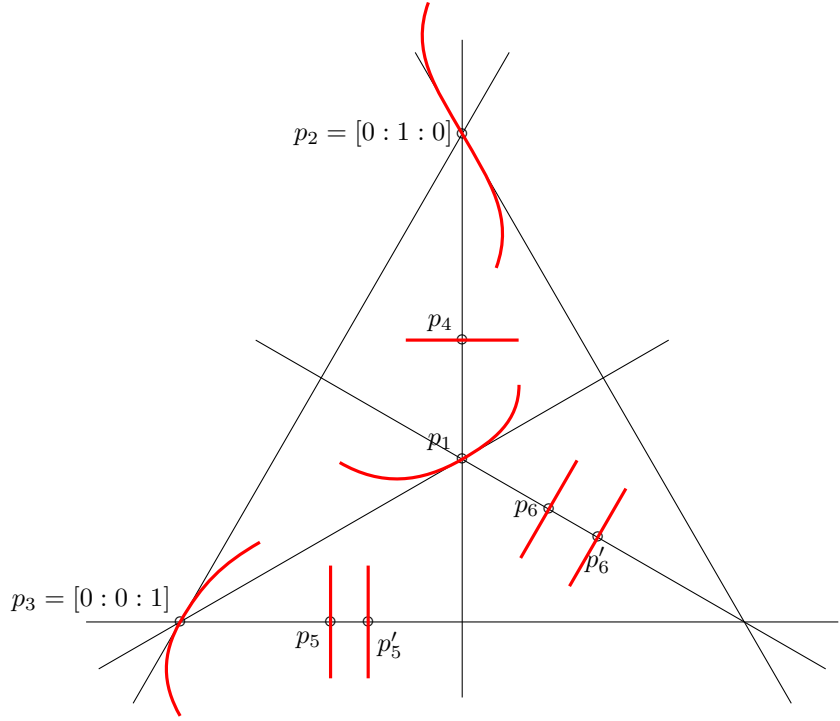
\begin{figure}
      \centering
      \begin{tikzpicture}[scale=2.5]
      \draw (-120:.5)--++(60:4);
      \draw (60:3)++(120:.5)--++(-60:4);
      \draw (0:-.5)--++(0:4);
      \draw (-150:.5)--++(30:3.5);
      \draw (60:3)++(90:.5)--++(-90:3.5);
      \draw (0:3)++(-30:.5)--++(150:3.5);
      \draw (0,0) node {$\circ$} node[above left] {$p_3=[0:0:1]$};
      \draw[red,very thick] (45:.6) to[out=-150,in=60] (0,0) to[out=-120,in=120] (-90:0.5);
      \draw (30:1.732) node {$\circ$} node[above left] {$p_1$};
      \draw[red,very thick] (45:1.2) to[out=-30,in=-150] (30:1.732) to[out=30,in=-90] (35:2.2);
      \draw (60:3) node {$\circ$} node[left] {$p_2=[0:1:0]$};
      \draw[red,very thick,xshift=3cm,xscale=-1] (63:3.7) to[out=-70,in=60] (60:3) to[out=-120,in=110] (55:2.3);
      \draw (1.5,1.5) node {$\circ$} node[above left] {$p_4$};
      \draw[red,very thick] (1.2,1.5)--(1.8,1.5);
      \draw (0.8,0) node {$\circ$} node[below left] {$p_5$};
      \draw[red,very thick] (0.8,-0.3)--(0.8,0.3);
      \draw (1,0) node {$\circ$} node[below right] {$p_5'$};
      \draw[red,very thick] (1,-0.3)--(1,0.3);
      \draw (3,0)++(150:1.2) node {$\circ$} node[left] {$p_6$};
      \draw[red,very thick] (3,0)++(150:1.2)--++(60:0.3);
      \draw[red,very thick] (3,0)++(150:1.2)--++(-120:0.3);
      \draw (3,0)++(150:0.9) node {$\circ$} node[below] {$p_6'$};
      \draw[red,very thick] (3,0)++(150:0.9)--++(60:0.3);
      \draw[red,very thick] (3,0)++(150:0.9)--++(-120:0.3);
    \end{tikzpicture}
      \caption{The result of blowing down the four blue curves in Figure \ref{fig:PictureOfCurve}}
      \label{fig:BlownDownPictureOfCurve}
  \end{figure}

  We may choose coordinates on $\P^2$ such that the $\psi_3(124|35)$ is identified with $[0:0:1],$ $\psi_3(145|23)$ is identified with $[0:1:0],$ $\psi_3(13|245)$ is identified with $[0:0:1],$
  and $\psi_3(125|34)$ is identified with $[1:1:1].$ Taking into account the new tangency conditions, the resulting picture of $\P^2$ is shown in Figure \ref{fig:BlownDownPictureOfCurve}. In this picture, we have coordinates $x/z=\CR(3,4,5,1)$ and $y/z=\CR(5,2,3,4)$, and images under $\psi_3\circ\pi_1$:
  \begin{align*}
      p_1&\mapsto[1:1:1]&p_2&\mapsto[0:1:0]&p_3&\mapsto[0:0:1]\\
      p_4&\mapsto[1:2:1]&p_5,p_5'&\mapsto[\mp i:0:1]&p_6,p_6'&\mapsto[\textstyle\frac{-1\mp\sqrt{5}}{2}:1:1]\\
      p_7,p_7'&\mapsto[\textstyle\frac{1\mp\sqrt{5}}{2}:\frac{3\mp\sqrt{5}}{2}:1].
  \end{align*}
  
  Let $\cE=\psi_3(\pi_1(\overline{\Per_{2,5}})).$ By B\'ezout's theorem, $\cE$ is a cubic curve. The locations of $p_1,\ldots,p_6'$ and tangency conditions give 12 linear equations in the 10 coefficients of a cubic curve. The equations are consistent, giving a confirmation of our computations up to this point, and the equation of $\cE$ is \begin{align*}
      x^3+y^2z-3xyz+xz^2=0.
  \end{align*}
  This equation is already in Weierstrass form, and is isomorphic to the elliptic curve 17a4 in the $L$-functions and Modular Forms Database via the change of coordinates $$[x:y:z]\mapsto[-x+z:-x+y+2z:z].$$ We note the $J$-invariant and periods of $\cE$:
  \begin{align*}
      J&=\frac{35937}{17}&\omega_1&\approx3.09416&\omega_2&\approx2.74574.
  \end{align*}
  
  Finally, we must show that $\pi_1|_{\overline{\Per_{2,5}}}$ is an embedding. To see this, note that $(\pi_1)^{-1}(\pi_1(p_2))=\{p_2\}$ (for if there were other points in the preimage, they would have shown up in Section \ref{sec:C5InHBar}). Furthermore, $p_2$ is not a ramification point of $\pi_1|_{\overline{\Per_{2,5}}}$, as follows. At $p_2,$ $s_1\sim\CR(C,a_1,a_{2,0},a_{4,0},a_{3,0})$ is a local coordinate along $\overline{\Per_{2,5}}$. By definition, $\pi_1$ preserves $s_1$, so the derivative of $\pi_1|_{\overline{\Per_{2,5}}}$ is nonvanishing at $p_2.$ It follows that $\pi_1|_{\overline{\Per_{2,5}}}$ is injective over $\pi_1(p_2)$, hence generically injective, hence --- since $\overline{\Per_{2,5}}$ and $\cE$ are smooth curves --- an embedding. As $\overline{\Per_{2,5}}$ is transverse to the four exceptional curves of $\psi_3,$ we have
  \begin{thm}
  $\overline{\Per_{2,5}}$ is isomorphic (as a $\Q$-variety) to the smooth plane curve $Z(x^3+y^2z-3xyz+xz^2)\subseteq\P^2.$
  \end{thm}

 \begin{rem}
 We briefly investigate the group structure on $\overline{\Per_{2,5}}.$ To do so, one must make an arbitrary choice of identity element $O$. One can check that $\overline{\Per_{2,5}}(\Q)=\{p_1,p_2,p_3,p_4\}$, and choosing any of these as $O$ yields $\overline{\Per_{2,5}}(\Q)\cong\Z/4\Z$. Thus if it is decided that $O$ should be a $\Q$-point, then $\overline{\Per_{2,5}}/\overline{\Per_{2,5}}(\Q)$ is canonically a group.
 Here are some computations via Sage:
  \begin{itemize}
      \item In $\overline{\Per_{2,5}}(\Q(i))/\overline{\Per_{2,5}}(\Q)\cong\Z$, we have $p_5'=-p_5\ne0.$ We do not know if $p_5$ is a generator.
      \item In $\overline{\Per_{2,5}}(\Q(\sqrt{5}))/\overline{\Per_{2,5}}(\Q)\cong\Z$, we have $p_6=-p_6'=-p_7=p_7'\ne0.$ Again, we do not know if $p_5$ is a generator.
  \end{itemize}
 \end{rem}

\subsection{An analog of the Mandelbrot set in \texorpdfstring{$\overline{\Per_{2,5}}$}{Per\_\{2,5\}-bar}}\label{sec:Mandelbrot}
In this section we work over $\C$. Let $\cP_2$ denote the moduli space of quadratic polynomial maps $\C\to\C$ up to conjugation, and recall that $\cP_2\cong\C$. An element $f\in\cP_2$ should be thought of as a rational function $f:\P^1\to\P^1$ with a critical fixed point $z_1$. (There is one other critical point $z_2$.)

There are two equivalent definitions of the Mandelbrot set $\mathscr{M}$ in $\cP_2$:
\begin{enumerate}
    \item $\mathscr{M}=\{f:\text{The Julia set of $f$ is connected}\}.$
    \item $\mathscr{M}=\{f:\text{$z_1$ is not in the closure of $\{f^n(z_2):n\ge0\}$}\}.$
\end{enumerate}
These two definitions make sense (though they are no longer equivalent) for any moduli space of rational maps $f:\P^1\to\P^1$ with exactly two critical points --- in particular, for $\Per_{2,5}$. It follows from \cite{Milnor2000} that for \textit{every} point $f\in\Per_{2,5}$, the Julia set of $f$ is connected. However, we may still consider $$\mathscr{M}_{2,5}:=\{f\in\Per_{2,5}:\text{$z_1$ is not in the closure of $\{f^n(z_2):n\ge0\}$}\}\subseteq\overline{\Per_{2,5}}.$$ (That is, $\mathscr{M}_{2,5}$ consists of rational maps $f\in\Per_{2,5}$ such that $z_2$ is not attracted to the specified 5-periodic orbit.) Since $\overline{\Per_{2,5}}\supseteq\Per_{2,5}$ is a smooth cubic curve, one may apply the inverse Weierstrass $\wp$-function to tile the complex plane $\C$ with copies of $\overline{\Per_{2,5}}$. Figure \ref{fig:Mandelbrot} on page \pageref{fig:Mandelbrot} shows a computer-generated image of four copies of $\mathscr{M}_5$ in $\C$. Punctures are marked in red and labeled according to the classification in Section \ref{sec:LocalAnalysis}. PCF points of types I, II, III, and IV are drawn as small blue dots in the centers of certain colored regions. Note that there are infinitely many other PCF points; we have only drawn those where both critical points are in the same 5-cycle. These are precisely the points that get collapsed in pairs to nodes of $\Per_{2,5}$ in $\mathcal{M}_2.$

\section{Application 2: A genus formula for \texorpdfstring{$\Per_{d,4}$}{Per\_\{d,4\}}}\label{sec:GenusFormula}
In this Section we consider degree-$d$ bicritical maps with a 4-periodic critical point, letting $d$ be arbitrary. We prove:
\begin{thm}\label{thm:GenusFormula}
 $\Per_{d,4}$ is isomorphic to a smooth plane curve of degree $d$, punctured at $d^2$ points. In particular, $\Per_{d,4}$ has genus $\frac{(d-1)(d-2)}{2}$ and gonality $d-1$.
\end{thm}
 We will realize $\Per_{d,4}$ as a smooth degree-$d$ plane curve; the degree-genus formula then gives the genus above. The analysis is similar to that in Section \ref{sec:C5InHBar}; we give an abbreviated computation, emphasizing the aspects that differ from the case of $\Per_{2,5}$.
 
 \subsection{Enumeration and analysis of boundary strata}\label{sec:LocalAnalysis2} In parallel to Section \ref{sec:StrataEnumeration}, we may enumerate (by hand or on a computer) the boundary strata in $\Hbar_{d,4}$ that intersect $\tilde{\Delta}_{d,4}$. These are given in Table \ref{fig:TableOfBoundaryStrataBicritical}. Each map of trees shown may correspond to a list of multiple boundary strata; a single boundary stratum is given by specifying the roots of unity by which certain cross-ratios differ. For example, along each stratum of type $\gamma_1$, the two cross-ratios $\CR(a_*,a_1,a_{2,0},a_{3,0})$ and $\CR(a_*,a_1,a_{2,0},a_{4,0})$ are constant, and are distinct $d$-th roots of unity (neither of which is 1). The $(d-1)(d-2)$ distinct choices of these roots give rise to $(d-1)(d-2)$ boundary strata. Note that strata of type $\gamma_1$ appear only when $d>2$. One might therefore worry that the table is missing some strata that appear only for very large $d$. There is, however, an easy argument that this does not occur, due essentially to the fact that there are very few marked points.
 
 \begin{table}\centering
  
  \begin{tabular}{cccc}
    \begin{tikzpicture}
    \draw (-1.5,0) node {$\tau$};
\draw (-1.5,2) node {$\sigma$};
\draw
      (0,0)--(1,0);
\draw (0,0) node {$\bullet$};
\draw
      (0,0)--++(160:.4);
\draw (0,0)++(160:.4) node[left] {2};
\draw
      (0,0)--++(200:.4);
\draw (0,0)++(200:.4) node[left] {$*$};
\draw
      (1,0) node {$\bullet$};
\draw (1,0)--++(-40:.4);
\draw
      (1,0)++(-40:.4) node[right] {1};
\draw (1,0)--++(0:.4);
\draw
      (1,0)++(0:.4) node[right] {3};
\draw (1,0)--++(40:.4);
\draw
      (1,0)++(40:.4) node[right] {4};
\draw[->] (.5,1)--(.5,0.5);
\foreach \j in {0,1,...,5} {
      \draw (0,2)--(1,1.5+\j/5);
      \draw (1,1.5+\j/5) node {$\bullet$};
      };
\draw (1,2.5) node {$\bullet$};
\draw (1,2.5)--++(40:.4);
\draw (1,2.5)++(40:.4) node[right] {3};
\draw (1,2.3)--++(0:.4);
\draw (1,2.3)++(0:.4) node[right] {2};
\draw (1,2.1)--++(-40:.4);
\draw (1,2.1)++(-40:.4) node[right] {4};
      \draw (0,2) node {$\bullet$};
\draw[very thick]
      (0,2)--++(160:.4);
\draw (0,2)++(160:.4) node[left] {1};
      \draw[very thick] (0,2)--++(200:.4);
\draw (0,2)++(200:.4)
      node[left] {$*$};
\draw (1,1.5) node {$\bullet$};
    \end{tikzpicture}&
    \begin{tikzpicture}
\draw
      (0,0)--(1,0);
\draw (0,0) node {$\bullet$};
\draw
      (0,0)--++(160:.4);
\draw (0,0)++(160:.4) node[left] {2};
\draw
      (0,0)--++(200:.4);
\draw (0,0)++(200:.4) node[left] {$*$};
\draw
      (1,0) node {$\bullet$};
\draw (1,0)--++(-40:.4);
\draw
      (1,0)++(-40:.4) node[right] {1};
\draw (1,0)--++(0:.4);
\draw
      (1,0)++(0:.4) node[right] {3};
\draw (1,0)--++(40:.4);
\draw
      (1,0)++(40:.4) node[right] {4};
\draw[->] (.5,1)--(.5,0.5);
\foreach \j in {0,1,...,5} {
      \draw (0,2)--(1,1.5+\j/5);
      \draw (1,1.5+\j/5) node {$\bullet$};
      };
\draw (1,2.5) node {$\bullet$};
\draw (1,2.5)--++(40:.4);
\draw (1,2.5)++(40:.4) node[right] {3};
\draw (1,2.5)--++(0:.4);
\draw (1,2.5)++(0:.4) node[right] {2};
\draw (1,2.5)--++(-40:.4);
\draw (1,2.5)++(-40:.4) node[right] {4};
      \draw (0,2) node {$\bullet$};
\draw[very thick]
      (0,2)--++(160:.4);
\draw (0,2)++(160:.4) node[left] {1};
      \draw[very thick] (0,2)--++(200:.4);
\draw (0,2)++(200:.4)
      node[left] {$*$};
\draw (1,1.5) node {$\bullet$};
    \end{tikzpicture}&
    \begin{tikzpicture}
\draw
      (0,0)--(1,0);
\draw (0,0) node {$\bullet$};
\draw
      (0,0)--++(140:.4);
\draw (0,0)++(140:.4) node[left] {2};
\draw
      (0,0)--++(180:.4);
\draw (0,0)++(180:.4) node[left] {4};
\draw
      (0,0)--++(220:.4);
\draw (0,0)++(220:.4) node[left] {$*$};
\draw
      (1,0) node {$\bullet$};
\draw (1,0)--++(-20:.4);
\draw
      (1,0)++(-20:.4) node[right] {1};
\draw (1,0)--++(20:.4);
\draw
      (1,0)++(20:.4) node[right] {3};
\draw[->] (.5,1)--(.5,0.5);
\foreach \j in {0,1,...,5} {
      \draw (0,2)--(1,1.5+\j/5);
      \draw (1,1.5+\j/5) node {$\bullet$};
      };
\draw (1,2.5) node {$\bullet$};
\draw (1,2.5)--++(-20:.4);
\draw (1,2.5)++(-20:.4) node[right] {4};
\draw (1,2.5)--++(20:.4);
\draw (1,2.5)++(20:.4) node[right] {2};
      \draw (0,2) node {$\bullet$};
\draw[very thick]
      (0,2)--++(140:.4);
\draw (0,2)++(140:.4) node[left] {1};
\draw[very thick]
      (0,2)--++(180:.4);
\draw (0,2)++(180:.4) node[left] {3};
      \draw[very thick] (0,2)--++(220:.4);
\draw (0,2)++(220:.4)
      node[left] {$*$};
\draw (1,1.5) node {$\bullet$};
    \end{tikzpicture}&
    \begin{tikzpicture}
\draw
      (0,0)--(2,0);
\draw (0,0) node {$\bullet$};
\draw
      (0,0)--++(160:.4);
\draw (0,0)++(160:.4) node[left] {2};
\draw
      (0,0)--++(200:.4);
\draw (0,0)++(200:.4) node[left] {$*$};
\draw
      (1,0) node {$\bullet$};
\draw (1,0)--++(90:.4);
\draw
      (1,0)++(90:.4) node[right] {1};
      \draw
      (2,0) node {$\bullet$};
\draw (2,0)--++(-20:.4);
\draw
      (2,0)++(-20:.4) node[right] {3};
\draw (2,0)--++(20:.4);
\draw
      (2,0)++(20:.4) node[right] {4};
\draw[->] (.5,1)--(.5,0.5);
\foreach \j in {0,1,...,5} {
      \draw (0,2)--(1,1.5+\j/5);
      \draw (1,1.5+\j/5) node {$\bullet$};
      \draw (1,1.5+\j/5)--(2,1.5+\j/5);
      \draw (2,1.5+\j/5) node {$\bullet$};
      };
\draw (1,2.5) node {$\bullet$};
\draw (2,2.5)--++(20:.4);
\draw (2,2.5)++(20:.4) node[right] {3};
\draw (2,2.3)--++(-20:.4);
\draw (2,2.3)++(-20:.4) node[right] {2};
\draw (1,2.5)--++(90:.4);
\draw (1,2.5)++(90:.4) node[above] {4};
      \draw (0,2) node {$\bullet$};
\draw[very thick]
      (0,2)--++(160:.4);
\draw (0,2)++(160:.4) node[left] {1};
      \draw[very thick] (0,2)--++(200:.4);
\draw (0,2)++(200:.4)
      node[left] {$*$};
\draw (1,1.5) node {$\bullet$};
    \end{tikzpicture}
    \\
    $(d-1)(d-2)$ strata&1 stratum&1 stratum&$d-1$ strata\\
    1 puncture each&$1$ puncture&$(d-1)$ punctures&1 puncture each\\
    $(\gamma_1)$&$(\gamma_2)$&$(\gamma_3)$&$(\gamma_4)$\\
    &&&\\\hline
    \begin{tikzpicture}
    \draw (-1.5,0) node {$\tau$};
\draw (-1.5,2) node {$\sigma$};
\draw
      (0,0)--(2,0);
\draw (0,0) node {$\bullet$};
\draw
      (0,0)--++(160:.4);
\draw (0,0)++(160:.4) node[left] {2};
\draw
      (0,0)--++(200:.4);
\draw (0,0)++(200:.4) node[left] {$*$};
\draw
      (1,0) node {$\bullet$};
\draw (1,0)--++(90:.4);
\draw
      (1,0)++(90:.4) node[right] {3};
      \draw
      (2,0) node {$\bullet$};
\draw (2,0)--++(-20:.4);
\draw
      (2,0)++(-20:.4) node[right] {1};
\draw (2,0)--++(20:.4);
\draw
      (2,0)++(20:.4) node[right] {4};
\draw[->] (.5,1)--(.5,0.5);
\foreach \j in {0,1,...,5} {
      \draw (0,2)--(1,1.5+\j/5);
      \draw (1,1.5+\j/5) node {$\bullet$};
      \draw (1,1.5+\j/5)--(2,1.5+\j/5);
      \draw (2,1.5+\j/5) node {$\bullet$};
      };
\draw (1,2.5) node {$\bullet$};
\draw (2,2.5)--++(20:.4);
\draw (2,2.5)++(20:.4) node[right] {3};
\draw (2,2.3)--++(-20:.4);
\draw (2,2.3)++(-20:.4) node[right] {4};
\draw (1,2.5)--++(90:.4);
\draw (1,2.5)++(90:.4) node[above] {2};
      \draw (0,2) node {$\bullet$};
\draw[very thick]
      (0,2)--++(160:.4);
\draw (0,2)++(160:.4) node[left] {1};
      \draw[very thick] (0,2)--++(200:.4);
\draw (0,2)++(200:.4)
      node[left] {$*$};
\draw (1,1.5) node {$\bullet$};
    \end{tikzpicture}
    &
     \begin{tikzpicture}
       \draw (0,0)--(1,0);
\draw (0,0) node {$\bullet$};
\draw
       (0,0)--++(160:.4);
\draw (0,0)++(160:.4) node[left] {3};
\draw
       (0,0)--++(200:.4);
\draw (0,0)++(200:.4) node[left] {$*$};
       \draw (1,0) node {$\bullet$};
\draw (1,0)--++(-40:.4);
\draw
       (1,0)++(-40:.4) node[right] {1};
\draw (1,0)--++(0:.4);
\draw
       (1,0)++(0:.4) node[right] {2};
\draw (1,0)--++(40:.4);
\draw
       (1,0)++(40:.4) node[right] {4};
\draw[->] (.5,1)--(.5,0.5);
       \draw[very thick] (1,1.5)--(0,1.5);
\draw (0,1.5) node {$\bullet$};
\draw
       (0,1.5)--++(160:.4);
\draw (0,1.5)++(160:.4) node[left] {2};
       \draw[very thick] (0,1.5)--++(200:.4);
\draw (0,1.5)++(200:.4)
       node[left] {$*$};
\draw (1,1.5) node {$\bullet$};
\draw
       (1,1.5)--++(-40:.4);
\draw (1,1.5)++(-40:.4) node[right] {4};
       \draw[very thick] (1,1.5)--++(0:.4);
\draw (1,1.5)++(0:.4)
       node[right] {1};
\draw (1,1.5)--++(40:.4);
\draw
       (1,1.5)++(40:.4) node[right] {3};
     \end{tikzpicture}
    &
     \begin{tikzpicture}
\draw (0,0)--(1,0);
\draw (0,0) node
      {$\bullet$};
\draw (0,0)--++(160:.4);
\draw (0,0)++(160:.4)
      node[left] {4};
\draw (0,0)--++(200:.4);
\draw (0,0)++(200:.4)
      node[left] {$*$};
\draw (1,0) node
      {$\bullet$};
\draw (1,0)--++(-40:.4);
\draw (1,0)++(-40:.4)
      node[right] {1};
\draw (1,0)--++(0:.4);
\draw (1,0)++(0:.4)
      node[right] {2};
\draw (1,0)--++(40:.4);
\draw (1,0)++(40:.4)
      node[right] {3};
\draw[->] (.5,1)--(.5,0.5);
\draw[very thick]
      (1,1.5)--(0,1.5);
\draw (0,1.5) node
      {$\bullet$};
\draw (0,1.5)--++(160:.4);
\draw (0,1.5)++(160:.4)
      node[left] {3};
\draw[very thick] (0,1.5)--++(200:.4);
\draw
      (0,1.5)++(200:.4) node[left]
      {$*$};
\draw (1,1.5) node
      {$\bullet$};
\draw (1,1.5)--++(-40:.4);
\draw (1,1.5)++(-40:.4)
      node[right] {4};
\draw[very thick] (1,1.5)--++(0:.4);
\draw
      (1,1.5)++(0:.4) node[right] {1};
\draw (1,1.5)--++(40:.4);
      \draw (1,1.5)++(40:.4) node[right] {2};
    \end{tikzpicture}&
                 \begin{tikzpicture}
                         \draw (0,0)--(1,0);
\draw (0,0) node
                         {$\bullet$};
\draw (0,0)--++(160:.4);
\draw
                         (0,0)++(160:.4) node[left] {1};
\draw
                         (0,0)--++(200:.4);
\draw (0,0)++(200:.4)
                         node[left] {$*$};
\draw (1,0) node
                         {$\bullet$};
\draw (1,0)--++(-40:.4);
\draw
                         (1,0)++(-40:.4) node[right] {2};
\draw
                         (1,0)--++(0:.4);
\draw (1,0)++(0:.4)
                         node[right] {3};
\draw (1,0)--++(40:.4);
                         \draw (1,0)++(40:.4) node[right] {4};
\draw[->] (.5,1)--(.5,0.5);
                         \draw[very thick] (1,1.5)--(0,1.5);
\draw (0,1.5) node
                         {$\bullet$};
\draw (0,1.5)--++(160:.4);
\draw
                         (0,1.5)++(160:.4) node[left] {4};
\draw[very
                         thick] (0,1.5)--++(200:.4);
\draw
                         (0,1.5)++(200:.4) node[left] {$*$};
\draw
                         (1,1.5) node {$\bullet$};
\draw[very thick]
                         (1,1.5)--++(-40:.4);
\draw (1,1.5)++(-40:.4)
                         node[right] {1};
\draw (1,1.5)--++(0:.4);
                         \draw (1,1.5)++(0:.4) node[right] {2};
                         \draw (1,1.5)--++(40:.4);
\draw
                         (1,1.5)++(40:.4) node[right] {3};
                       \end{tikzpicture}
                       \\
    $(d-1)$ strata&1 stratum&1 stratum&1 stratum\\
    1 puncture each&0 punctures&0 punctures&0 punctures\\
    $(\gamma_5)$&$(\gamma_{\mathrm{I}})$&$(\gamma_{\mathrm{II}})$&$(\gamma_{\mathrm{III}})$\\
    &&&\\
  \end{tabular}
  
  \caption{Strata that intersect $\overline{\Per_{d,4}}$. The case $d=6$ is shown, but it is an easy argument that the same strata types appear for any $d\ge3$.}
  \label{fig:TableOfBoundaryStrataBicritical}  
\end{table}

These boundary strata may be analyzed as in Section \ref{sec:LocalAnalysis}; we give only a sample computation, as the others may be easily checked.

\medskip

\noindent\textit{Analysis of $\gamma_3$.}
The combinatorial type $\gamma_3$ corresponds to a 1-dimensional locally-closed stratum $R_3$.
The two functions $s_1=CR(a_1,a_{2,0},a_{3,0},a_{4,0})$ and
$s_2=CR(a_*,a_1,a_{2,0},a_{3,0})$ give \'etale coordinates on an \'etale chart
$U_3$ containing $R_3$. Note that $s_1$ is a
node-smoothing parameter and $s_2$ is a cross-ratio parameter on the
left vertex of $\sigma,$ so $R_3$ is defined by $s_1$ on $U_3$. (Note also that on $R_3$, we have $s_2^d\ne1$.) The
map $(s_1,s_2):\Hbar_{d,4}\to\A^2$ is an embedding, i.e. $U_3$ is
actually a Zariski chart.

Locally, $\tilde{\Delta}_{d,4}$ is defined by
$h:=CR(a_1,a_{2,0},a_{3,0},a_{4,0})-CR(b_1,b_2,b_3,b_4).$ The cross-ratio
relations imply the following expansion:
\begin{align}
  h&=s_1\left(\frac{(d+1)s_2^d-ds_2^{d-1}-1}{s_2^d-1}\right)+O(s_1^2).
\end{align}
As noted above, the case $s_2^d-1=0$ corresponds to a degeneration of $R_3$, so this is indeed a regular function on $R_3$.

We conclude that the intersection of $\mathcal{C}_{d,4}$ with $R_3$ is cut
out by ($s_1=0$ and) $$g:=(d+1)s_2^d-ds_2^{d-1}-1=0.$$ Note that
$s_2=1$ is a root of $g$; note also that all other $d$-th roots of
unity are \emph{not} roots of $g$.

The discriminant of $g$ is given by $d^d((d+1)^{d-1}+(d-1)^{d+1});$
this follows directly from the form of the Sylvester matrix of $g$ and
$\pderiv{g}{s_2}$. As the discriminant is clearly nonzero, $g$ has $d$
distinct roots (including $s_2=1$, corresponding to a point not in
$R_3$), so we conclude that $\mathcal{C}_{d,4}$ intersects $R_3$ in
exactly $d-1$ distinct reduced points; as usual, these must be smooth points of $\mathcal{C}_{d,4}$.

\medskip

For our purposes, it will be useful to consider the image of $\mathcal{C}_{d,4}$ under the map $\pi_2:\Hbar_{d,4}\to\Mbar_{0,5}$ that remembers the \emph{source} curve. From Section \ref{sec:HnBarCoords}, $\pi_2$ is birational, and is an isomorphism away from certain boundary strata where the source curve is unstable. (An example where this happens is $\gamma_1$, which corresponds to a collection on 1-dimensional boundary strata, each of which collapses to a single point in the interior of $\Mbar_{0,5}.$) In particular, $\pi_2$ is injective in a neighborhood of $R_3$, so the $d-1$ smooth points of $\mathcal{C}_{d,4}$ above map to $d-1$ distinct smooth points of $\pi_2(\mathcal{C}_{d,4})$, all of which lie on the boundary stratum $13{*}|24\subseteq\Mbar_{0,5}$.

\subsection{Summary of boundary strata analyses.}  Essentially identical arguments for the rest of the boundary strata give the following conclusions:
\begin{itemize}
    \item Combinatorial type $\gamma_{1}$ corresponds to $(d-1)(d-2)$ distinct smooth points of $\mathcal{C}_{d,4}$ --- one for each stratum. These map to $(d-1)(d-2)$ distinct smooth points of $\pi_2(\mathcal{C}_{d,4})$ in the interior of $\Mbar_{0,5}$, where the cross-ratios $\CR(a_1,a_*,a_{2,0},a_{3,0})$ and $\CR(a_1,a_*,a_{2,0},a_{4,0})$ are determined by the stratum.
    \item Combinatorial type $\gamma_{2}$ corresponds to a single smooth point of $\mathcal{C}_{d,4}$, which maps to a smooth point of  $\pi_2(\mathcal{C}_{d,4})$ on $1{*}|234$.
    \item Combinatorial type $\gamma_{3}$ corresponds to $d-1$ distinct smooth points of $\mathcal{C}_{d,4}$, which map to $d-1$ smooth points of $\pi_2(\mathcal{C}_{d,4})$ on $13{*}|24$.
    \item Combinatorial type $\gamma_{4}$ corresponds to $d-1$ distinct smooth points of $\mathcal{C}_{d,4}$ --- one for each stratum. These map to $d-1$ smooth points of $\pi_2(\mathcal{C}_{d,4})$ on $12{*}|34$, where the cross-ratios $\CR(a_1,a_*,a_{2,0},a_{3,0})$ are distinct $d$th roots of unity not equal to 1.
    \item Combinatorial type $\gamma_{5}$ corresponds to $d-1$ distinct smooth points of $\mathcal{C}_{d,4}$ --- one for each stratum. These map to $d-1$ smooth points of $\pi_2(\mathcal{C}_{d,4})$ on $14{*}|23$, where the cross-ratios $\CR(a_1,a_*,a_{4,0},a_{2,0})$ are distinct $d$th roots of unity not equal to 1.
    \item Combinatorial type $\gamma_{\mathrm{I}}$ corresponds to $d$ distinct smooth points of $\mathcal{C}_{d,4}$, which map to $d$ smooth points of $\pi_2(\mathcal{C}_{d,4})$ on $134|2{*}$.
    \item Combinatorial type $\gamma_{\mathrm{II}}$ corresponds to $d$ distinct smooth points of $\mathcal{C}_{d,4}$, which map to $d$ smooth points of $\pi_2(\mathcal{C}_{d,4})$ on $124|3{*}$.
    \item Combinatorial type $\gamma_{\mathrm{III}}$ corresponds to $d$ distinct smooth points of $\mathcal{C}_{d,4}$, which map to $d$ smooth points of $\pi_2(\mathcal{C}_{d,4})$ on $123|4{*}$.
\end{itemize}
In particular, we may conclude that $\pi_2(\mathcal{C}_{d,4})$ is smooth.
 
\begin{figure}
    \centering
    \begin{tikzpicture}[scale=0.9]
      \draw (-30:.5) to[out=-90,in=45] (-90:3) to[out=-135,in=0, edge node={node[above] {$125|34$}}] (-135:4);
      \draw (90:.5) to[out=30,in=165] (30:3) to[out=-15,in=120, edge node={node[below left=.1cm] {$13{*}|24$}}] (-15:4);
      \draw (210:.5) to[out=150,in=285] (150:3) to[out=105,in=240, edge node={node[below right] {$14{*}|23$}}] (105:4);
      \draw[thick,blue] (-125:4) to[out=90,in=0] (-160:4) to[out=180,in=-90, edge node={node[above left] {$12|34{*}$}}] (110:6);
      \draw[thick,blue] (-5:4) to[out=210,in=120] (-40:4) to[out=300,in=30, edge node={node[below] {$13|24{*}$}}] (230:6);
      \draw[thick,blue] (115:4) to[out=330,in=240] (80:4) to[out=420,in=150, edge node={node[right] {$14|23{*}$}}] (350:6);
      \draw (-210:2) to[out=-180,in=90] (-160:4) to[out=-90,in=90, edge node={node[left] {$123|4{*}$}}] (-120:6);
      \draw (-90:2) to[out=-60,in=210] (-40:4) to[out=30,in=210, edge node={node[below right] {$134|2{*}$}}] (0:6);
      \draw (30:2) to[out=60,in=330] (80:4) to[out=150,in=330, edge node={node[above right] {$124|3{*}$}}] (120:6);
      \draw[thick,blue] (0:1) to[out=90,in=0] (90:1) node[above=4pt,left=-3pt] {$1{*}|234$} to[out=180,in=90] (180:1) to[out=270,in=180] (270:1) to[out=0,in=270] (0:1);
      \foreach \i in {0,1,...,4} {
      \draw (-68-2*\i:1.2+\i/5) node {$\circ$} node[below=2pt] {};
      \draw[red,very thick] (-68-2*\i:1.2+\i/5)--++(0:.25);
      \draw[red,very thick] (-68-2*\i:1.2+\i/5)--++(180:.25);};
      \draw (-45:1.7) node {\textcolor{red}{$d-1$}};
      \foreach \i in {0,1,...,4} {
      \draw (152-2*\i:2.9+\i/10) node {$\circ$} node[below=2pt] {};
      \draw[red,very thick] (152-2*\i:2.9+\i/10)--++(5-\i/2:.5);
      \draw[red,very thick] (152-2*\i:2.9+\i/10)--++(185-\i/2:.5);};
      \draw (134:2.2) node {\textcolor{red}{$d-1$}};
      \foreach \i in {0,1,...,5} {
      \draw (180+2*\i:3.4+\i/15) node {$\circ$} node[below=2pt] {};
      \draw[red,very thick] (180+2*\i:3.4+\i/15)--++(-5-\i/2:.5);
      \draw[red,very thick] (180+2*\i:3.4+\i/15)--++(-185-\i/2:.5);};
      \foreach \i in {0,1,...,5} {
      \draw (45+2*\i:5.4+\i/15) node {$\circ$} node[below=2pt] {};
      \draw[red,very thick] (45+2*\i:5.4+\i/15)--++(0-\i/2:.5);
      \draw[red,very thick] (45+2*\i:5.4+\i/15)--++(180-\i/2:.5);};
      \foreach \i in {0,1,...,5} {
      \draw (70-2*\i:3.7-\i/15) node {$\circ$} node[below=2pt] {};
      \draw[red,very thick] (70-2*\i:3.7-\i/15)--++(30:.5);
      \draw[red,very thick] (70-2*\i:3.7-\i/15)--++(210:.5);};
      \draw (75:3) node {\textcolor{red}{$d$}};
      \draw (-32:5.2) node {\textcolor{red}{$d$}};
      \draw (188:2.8) node {\textcolor{red}{$d$}};
      \draw (60:6) node {\textcolor{red}{$(d-1)(d-2)$}};
      \draw (-120:1) node {$\circ$} node[below=2pt,left=0pt] {};
      \draw[red,very thick] (-120:.5)--(-120:1.5);
      \draw (20:3.33) node {$\circ$} node[above=2pt,left=0pt] {};
      \draw[red,very thick] (20:3.33)--++(50:.5);
      \draw[red,very thick] (20:3.33)--++(230:.5);
      \draw (22.5:3.22) node {$\circ$} node[above=2pt,left=0pt] {};
      \draw[red,very thick] (22.5:3.22)--++(55:.5);
      \draw[red,very thick] (22.5:3.22)--++(235:.5);
      \draw (25:3.14) node {$\circ$} node[above=2pt,left=0pt] {};
      \draw[red,very thick] (25:3.14)--++(60:.5);
      \draw[red,very thick] (25:3.14)--++(240:.5);
      \draw (27.5:3.08) node {$\circ$} node[above=2pt,left=0pt] {};
      \draw[red,very thick] (27.5:3.08)--++(65:.5);
      \draw[red,very thick] (27.5:3.08)--++(245:.5);
      \draw (30:3) node {$\circ$} node[above=8pt,left=0pt] {};
      \draw[red,very thick] (30:3)--++(70:.5);
      \draw[red,very thick] (30:3)--++(250:.5);
      \draw (12:2.5) node {\textcolor{red}{$d-1$}};
      \foreach \i in {0,1,...,5} {
      \draw (-32+2*\i:4.21+\i/15) node {$\circ$} node[below=2pt] {};
      \draw[red,very thick] (-32+2*\i:4.21+\i/15)--++(130+4*\i:.5);
      \draw[red,very thick] (-32+2*\i:4.21+\i/15)--++(-50+4*\i:.5);};
      \draw (-32:5.2) node {\textcolor{red}{$d$}};
    \end{tikzpicture}
    
    \caption{Diagram of how $\pi_2(\overline{\Per}_{2,5})\subseteq\Mbar_{0,5}$ intersects boundary strata}
    \label{fig:PictureOfC4}
  \end{figure}
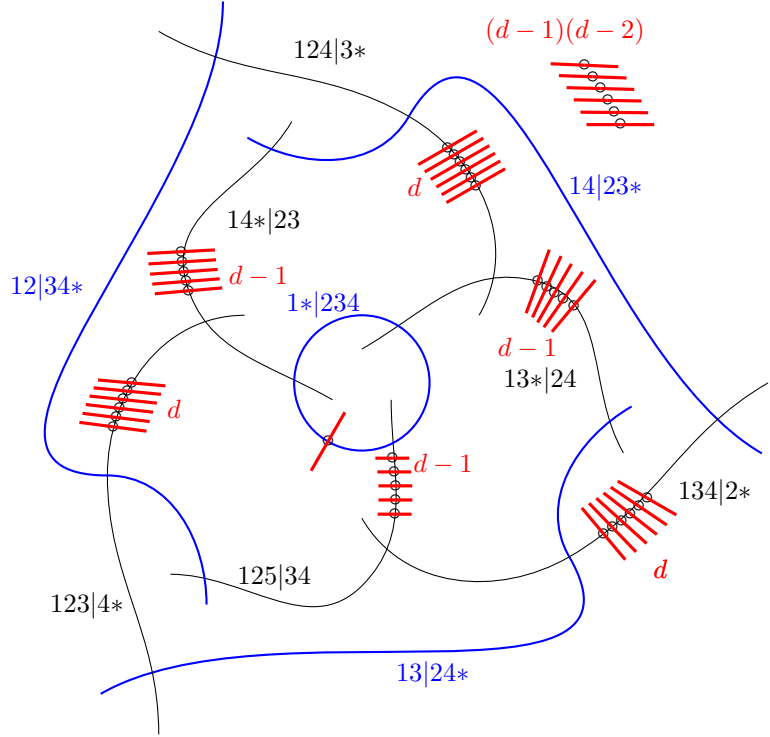 
  
  As in Section \ref{sec:C5InHBar}, the types $\gamma_{\mathrm{I}}$, $\gamma_{\mathrm{II}}$, and $\gamma_{\mathrm{III}}$ correspond to points (rather than punctures) of $\Per_{d,4}$, namely rational functions on $\P^1$ in which both critical points are in the same 4-cycle. Adding up all the punctures above, we get $d^2$ points in total, as claimed.

\subsection{Blowing down to \texorpdfstring{$\P^2$}{P2}}
To prove Theorem \ref{thm:GenusFormula}, we blow down $\Mbar_{0,5}$ to $\P^2$ along the four boundary divisors indicated in Figure \ref{fig:PictureOfC4}; this is the map $\psi_1$, analogous to Section \ref{sec:CubicCurve}. Observe that of these boundary divisors, only $1{*}|234$ intersects $\pi_2(\mathcal{C}_{d,4})$, and that in only one smooth point; it follows that $\psi_1(\pi_2(\mathcal{C}_{d,4}))\subseteq\P^2$ is smooth. Observe also that any noncontracted boundary divisor meets $\psi_1(\pi_2(\mathcal{C}_{d,4}))$ in exactly $d$ points; that these numbers agree is a check on the consistency of our calculations. (Note that e.g. $13{*}|24$ picks up the extra point to $1{*}|234$ to make up $d$ points altogether.) As $\psi_1\circ\pi_2$ is an isomorphism on $\mathcal{C}_{d,4}$, we conclude Theorem \ref{thm:GenusFormula}.

\begin{rem}
In particular, the gonality of $\mathcal{C}_{d,4}$, namely $d-1$, is far lower than the gonality of a general curve of genus $\frac{(d-1)(d-2)}{2}$, namely $\floor{\frac{d^2-3d+8}{2}}$. More generally, it is interesting to study how ``special'' COR curves are within their respective moduli spaces.
\end{rem}
 
\bibliographystyle{alpha} \bibliography{BigRefs.bib}

\end{document}